\documentclass[10pt,leqno]{amsart} 

   \usepackage[centertags]{amsmath}
   \usepackage{amsfonts}
   \usepackage{amsmath}
   \usepackage{amssymb}
   \usepackage{amsthm}
   \usepackage{newlfont}
   \usepackage[latin1]{inputenc}
\usepackage{xcolor}

\allowdisplaybreaks[3]

\theoremstyle{plain}
\newtheorem{Theo}{Theorem}[section]

\newtheorem{Cor}[Theo]{Corollary}
\newtheorem{Lem}[Theo]{Lemma}

\theoremstyle{definition}
\newtheorem{Def}{Definition}[section]

\theoremstyle{remark}
\newtheorem{Rem}{Remark}[section]

\newcommand{\lge}{\mathfrak{b}}

\newcommand{\lp}{\mathcal{L}\text{-}\mathcal{P}}

\newcommand{\RR}{\mathbb{R}}
\newcommand{\ZZ}{\mathbb{Z}}
\newcommand{\NN}{\mathbb{N}}
\newcommand{\CC}{\mathbb{C}}

\newcommand{\nr}{\operatorname{nr}}

\numberwithin{equation}{section}

\newcommand{\card}{\operatorname{card}}

\newcommand\Ll{{\mathfrak L}}
\newcommand\La{{\mathcal L}}

\newcommand\Be{\mathfrak{b}}

\newcommand\sign{\operatorname{sign}}

\newcommand*\pFq{\begingroup
        \dopFq
}
\catcode`\,12
\def\dopFq#1#2#3#4#5{%
        {}_{#1}F_{#2}\biggl(\genfrac..{0pt}{}{#3}{#4};#5\biggr)%
        \endgroup
}

\title{Zeros of linear combinations of Laguerre polynomials}
\author{Antonio J. Dur\'an}
\address{Departamento de An\'a\-li\-sis Mate\-m\'a\-ti\-co and IMUS,
        Universidad de Sevilla,
        41080 Sevilla, Spain}
\email{duran@us.es}
   \date{}

   \thanks{This research was partially supported by PID2021-124332NB-C21
(Minis\-te\-rio de Cien\-cia e Inno\-va\-ci\'on and Feder Funds (European Union)), and
FQM-262 (Jun\-ta de Anda\-lu\-c\'ia).}



\begin{document}
   \maketitle

\begin{abstract}
We study the number of real zeros of finite combinations of $K+1$ consecutive normalized Laguerre polynomials of the form
$$
q_n(x)=\sum_{j=0}^K\gamma_j\tilde L^\alpha_{n-j}(x),\quad n\ge K,
$$
where $\gamma_j$, $j=0,\cdots ,K$, are real numbers with $\gamma_0=1$, $\gamma_K\not =0$. We consider four different normalizations of Laguerre polynomials: the monic Laguerre polynomials $\hat L_n^\alpha$, the polynomials $\La_n^\alpha=n!L_n^\alpha/(1+\alpha)_n$ (so that $\La_n^\alpha(0)=1$), the standard Laguerre polynomials $(L_n^\alpha)_n$ and the Brenke normalization $L_n^\alpha/(1+\alpha)_n$.
We show the key role played by the polynomials $Q(x)=\sum_{j=0}^K(-1)^j\gamma_j(x)_{K-j}$ and $P(x)=\sum_{j=0}^K\gamma_jx^{K-j}$ to solve this problem: $Q$ in the first case and $P$ in the second, third and forth cases.
In particular, in the first case, if all the zeros of the polynomial $Q$ are real and less than $\alpha+1$, then all the zeros of $q_n$, $n\ge K$, are positive. In the other cases,
if all the zeros of $P$ are real then all the zeros of $q_n$, $n\ge K$, are also real. If $P$ has $m>1$ non-real zeros, there are important differences between the four cases. For instance
in the first case, $q_n$ has still only real zeros for $n$ big enough, but in the fourth case $q_n$ has exactly $m$ non-real zeros for $n$ big enough.
\end{abstract}

\section{Introduction}
Consider the following result recently proved in the paper \cite{Dur0}.

\begin{Theo}[Corollary 3.4 of \cite{Dur0}] For any positive measure $\mu$ in the real line (and infinitely many points in its support), there always exists a sequence of
orthogonal polynomials $(p_n)_n$ with respect to $\mu$ such that for any positive integer $K$ and any $K+1$ real numbers
$\gamma_j$, $j=0,\cdots ,K$, with $\gamma_0=1$, $\gamma_K\not =0$, the polynomial
\begin{equation}\label{dqnv}
q_n(x)=\sum_{j=0}^K\gamma_jp_{n-j}(x),\quad n\ge K,
\end{equation}
has only real zeros for $n$ big enough (depending on $K$ and the $\gamma_j$'s).
\end{Theo}

Shohat \cite{Sho} was probably the first to observe that the orthogonality of the sequence $(p_n)_n$ implies that $q_n$ has at least $n-K$ real zeros in the convex hull of the support of $\mu$ (using the usual proof that $p_n$ has its $n$ zeros in the convex hull of the support of $\mu$). Some other related results on zeros of linear combinations of the form (\ref{dqnv}) can be found in \cite{Peh1,Peh2,IsSa,IsNo,INS,BDr,Fis,DJM,CSW,KRZ}.

We have to notice that the problem of studying the zeros of finite linear combinations of orthogonal polynomials of the form (\ref{dqnv})
is strongly dependent on the normalization of the polynomials $(p_n)_n$.
We have also proved in \cite{Dur0} that our result applies to the monic Laguerre polynomials.

The purpose of this paper is to show that the quasi-spectral properties of Laguerre polynomials allow to prove many more interesting and somehow unexpected results on the zeros of finite linear combinations of four different normalizations of this classical family of orthogonal polynomials: the monic Laguerre polynomials $(\hat L_n^\alpha)_n$, the Laguerre polynomials $(L_n^\alpha)_n$ themselves and the normalizations
\begin{align}\label{lag3x} \La_n^\alpha(x)=\frac{n!}{(1+\alpha)_n}L_n^\alpha(x),\\ \label{lag2nv}
\Ll_n^\alpha(x)=\frac{1}{(1+\alpha)_n}L_n^\alpha(x).
\end{align}
For a study of the Hermite case see \cite{Durh}.

Each one of these four normalizations of the Laguerre polynomials enjoys quasi-spectral properties with respect to a certain differential operator. Namely, define
\begin{align}\label{1op}
\Lambda_\alpha &=x \frac{d}{dx}+(\alpha-x),\\\label{2op}
\Upsilon_\alpha&=\frac{1}{\alpha}\left(x\frac{d}{dx}+\alpha\right),\\\label{3op}
\Omega_\alpha&=-\frac{d}{dx}\left(\alpha +x \frac{d}{dx}\right).
\end{align}
The first order differential operator (\ref{1op}) is the backward shift operator with respect to which the monic Laguerre polynomials enjoy the following quasi-spectral property
\begin{equation}\label{opli}
\Lambda_\alpha (\hat L^\alpha_n)=-\hat L^{\alpha-1}_{n+1}.
\end{equation}
We have used the terminology \textit{quasi-spectral} because the polynomial $\hat L^\alpha_n$ is not properly an eigenfunction of $\Lambda_\alpha$.

The normalized Laguerre polynomials (\ref{lag3x}) enjoy the following quasi-spectral property with respect to the first order differential operator (\ref{2op})
\begin{equation}\label{2opli}
\Upsilon_\alpha (\La^\alpha_n)=\La^{\alpha-1}_{n}.
\end{equation}
The normalized Laguerre polynomials (\ref{lag2nv}) enjoy in turn the following quasi-spectral property with respect to the second order differential operator (\ref{3op})
\begin{equation}\label{3opli}
\Omega_\alpha (\Ll^\alpha_n)=\Ll^{\alpha}_{n-1}.
\end{equation}
Notice that, in all the cases, the \textit{quasi-eigenvalues} are constant equal to $-1$ in the first case, and $1$ in the other cases. This is essential to guarantee the properties of the zeros of the polynomials $q_n$ (\ref{dqnv}) we will prove below.

As for the standard normalization of the Laguerre polynomials we have
\begin{equation}\label{4opli}
(L_n^\alpha)'(x)=-L^{\alpha+1}_{n-1}(x),\quad L_n^{\alpha-1}(x)=L_n^\alpha(x)-L_{n-1}^\alpha(x).
\end{equation}
In a different way, each of these quasi-spectral properties allow to study in depth the properties of the zeros of finite linear combinations of each normalized Laguerre families (which will also require different asymptotic properties of the Laguerre polynomials).

The content of the paper is as follows. In Section \ref{ssec4}, we consider the monic Laguerre polynomials. The starting point is the following result proved in \cite{Dur0}.

\begin{Cor}[Remark 2 in Section 4.2 of \cite{Dur0}]\label{1lagi}
Let $\alpha>-1$. Then for any positive integer $K$ and any finite set of $K+1$ real numbers $\gamma_j$, $j=0,\dots, K$, with $\gamma_0=1$ and $\gamma_K\not=0$, the polynomial
\begin{equation}\label{qnmli}
q_n(x)=\sum_{j=0}^K\gamma_j\hat L^\alpha _{n-j}(x)
\end{equation}
has only real and simple zeros for $n\ge \max\{\frac{1}{8}( \frac{23}{22}46^{K-1}\max\{|\gamma_j|,2\le j\le K\}-2|\alpha|),\alpha+9,2K\}$, and they interlace the zeros of $L^\alpha_{n-1}$. Moreover, for $n$ big enough all the zeros are positive.
\end{Cor}

Using the spectral property (\ref{1op}), we have improved Corollary \ref{1lagi} as follows. Consider the polynomial
\begin{equation}\label{dPl}
Q(x)=\sum_{j=0}^K(-1)^j\gamma_{j}(x)_{K-j},
\end{equation}
where $(u)_s$, $u\in \RR$ and $s\in \NN$, denotes the Pochhammer symbol: $(u)_s=u(u+1)\cdots (u+s-1)$.

\begin{Theo}\label{ldlc} Let $\alpha$ and $K$ be a real number $\alpha>-1$ and a positive integer, respectively.
Given $K+1$ real numbers $\gamma_j$, $j=0,\dots, K$, with $\gamma_0=1$ and $\gamma_K\not=0$, consider the polynomial $Q$ defined in (\ref{dPl}) and denote by $\theta_j$, $1\le j\le m$, the real zeros of $Q$ satisfying $\alpha+1\le \theta_j$, arranged in increasing order.
\begin{enumerate}
\item If the polynomial $Q$ has only real zeros then all the zeros of the polynomial $q_n$ (\ref{qnmli}) are real and simple for $n\ge n_1$, and positive and simple for $n\ge n_0$, where
$n_l= \max\{K,\lfloor\theta_i-\alpha+K\rfloor, 1\le i\le m-l\}$ (where, as usual, $\lfloor x\rfloor$ denotes the floor function which gives the greatest integer $\le x$).
In particular, if all the zeros of $Q$ are less than $\alpha+1$ then $n_0=n_1=K$.
\item If $Q$ has non-real zeros, then there exists a positive integer $n_0$, depending only on the non-real zeros of $Q$ and the real zeros $\theta_j$, $1\le j\le m$,
such that all the zeros of the polynomial $q_n$  are positive and simple for $n\ge n_0$.
\end{enumerate}
In all the cases, for $n$ big enough the zeros of $q_{n+1}$ interlace the zeros of $q_n$.
\end{Theo}

In order to prove this result, using the differential operator (\ref{opli}), we have introduced some new families of polynomials (related to the Bell polynomials),
which are interesting in themselves, whose zeros have very nice properties. We have proved that the polynomials $(q_n)_n$ (\ref{qnmli}) are particular cases of these families.

In Section \ref{secl33}, we study finite combinations of the normalized Laguerre polynomials (\ref{lag3x}) of the form
\begin{equation}\label{qnml5}
q_n(x)=\sum_{j=0}^K \gamma_j \La^\alpha_{n-j}(x).
\end{equation}
Let us note that $\La_n^\alpha(0)=1$.
In this case, we have to consider the polynomial
\begin{equation}\label{dPli}
P(x)=\sum_{j=0}^K\gamma_jx^{K-j}.
\end{equation}
We will need to assume that
\begin{equation}\label{jol}
\mbox{$P(z)\not =0$ for $z\in \{z\in \CC\setminus \RR:\vert z\vert \le  1\}$.}
\end{equation}

The spectral property (\ref{2opli}) leads to the following theorem.

\begin{Theo}\label{hithis} Let $\alpha> -1$. Given a positive integer $K$ and $K+1$ real numbers $\gamma_j$, $0\le j\le K$, with $\gamma_0=1$, $\gamma_K\not =0$, assume that the polynomial
(\ref{dPli}) satisfies (\ref{jol}). We denote by $N^{\nr}$ the number of non-real zeros of $P$ and by $N^1$ the number of real zeros of $P$ greater than $1$.
If $P(1)\not =0$, we then have:
\begin{enumerate}
\item If $N^{\nr}=0$, then $q_n$ (\ref{qnmlii}) has $n$ real and simple zeros for $n\ge K$, of which $N^{1}$ are negative, and the zeros of $q_{n+1}$ interlace the zeros of $q_n$.
\item Otherwise, there exists $n_0$ such that for $n\ge n_0$ the polynomial $q_n$ has exactly $n-N^{\nr}$ real and simple zeros, of which $N^{1}$ are negative, and the real zeros of $q_{n+1}$ interlace the real zeros of $q_n$.
\end{enumerate}
If $P$ has a zero at $x=1$ of multiplicity $s$, then $q_n$ has a zero of multiplicity $s$ at $x=0$ and then the zeros of $q_n/x^s$ behave as above.
\end{Theo}

In Section \ref{secll}, we study finite combinations of the Laguerre polynomials of the form
\begin{equation}\label{qnmlii}
q_n(x)=\sum_{j=0}^K \gamma_j L^\alpha_{n-j}(x).
\end{equation}
Using the spectral properties (\ref{4opli}), we have also prove a quite complete description of the zeros of the polynomials $q_n$.

\begin{Theo}\label{hithi} Given a positive integer $K$ and $K+1$ real numbers $\gamma_j$, $0\le j\le K$, with $\gamma_0=1$, $\gamma_K\not =0$, assume that the polynomial
$P$ (\ref{dPli}) satisfies (\ref{jol}). We denote by $N^{\nr}$ the number of non-real zeros of $P$ and by $N^1$ the number of real zeros of $P$ greater than $1$. We then have for $\alpha\ge K-N^{\nr}-1$:
\begin{enumerate}
\item If $N^{\nr}=0$, then $q_n$ (\ref{qnmlii}) has $n$ real and simple zeros for $n\ge K$, and the zeros of $q_{n+1}$ interlace the zeros of $q_n$. If $N^{1}=0$ then all the zeros are positive. If $N^{1}>0$ then for $n$ big enough  $q_n$ has exactly  $N^{1}$ negative zeros.
\item Otherwise, for $n$ big enough $q_n$ has exactly $n-N^{\nr}$ real and simple zeros of which $N^{1}$ are negative, and the real zeros of $q_{n+1}$ interlace the real zeros of $q_n$.
\end{enumerate}
\end{Theo}

In Section \ref{secl3}, we study finite combinations of the normalized Laguerre polynomials (\ref{lag2nv}) of the form
\begin{equation}\label{qnml3}
q_n(x)=\sum_{j=0}^K \gamma_j \Ll^\alpha_{n-j}(x).
\end{equation}
With this normalization of the Laguerre family, the polynomials $q_n$ are Brenke polynomials, in the sense that they can be defined by the following generating function (see \cite{drh})
\begin{equation}\label{genl}
\sum_{n=0}^\infty q_n(x)z^n=e^{z}z^KP(1/z)\pFq{0}{1}{-}{\alpha+1}{-xz},
\end{equation}
where the polynomial $P$ is defined above (\ref{dPli}).

Using that $(q_n)_n$ are Brenke polynomials and the quasi-spectral property (\ref{3opli}), we prove the following theorem.

\begin{Theo}\label{ultt} Assume that the polynomial $P$ (\ref{dPli}) has $N^{\nr}$ non-real zeros and $N^+$ positive zeros. Then for $\alpha>-1$ we have:
\begin{enumerate}
\item  $N^{\nr}=0$ if and only if $q_n$ (\ref{qnml3}) has only real zeros for all $n\ge 0$. Moreover, the zeros of $q_n$ are simple, the zeros of $q_{n+1}$ interlace the zeros of $q_n$ and for $n$ big enough, $q_n$ has exactly $N^+$ negative zeros.
\item If $N^{\nr}>0$ then there exists $n_0$ such that $q_n$ has exactly $n-N^{\nr}$ real zeros for $n\ge n_0$ of which $N^+$ are negative. For $\alpha\ge 0$, the real zeros of $q_{n+1}$ interlace the real zeros of $q_n$.
\end{enumerate}
\end{Theo}

We remark that the proofs of Theorems \ref{ldlc}, \ref{hithis}, \ref{hithi} and \ref{ultt} will use different approaches.

\medskip

In Section \ref{apen}, we consider the case $n=K$. Iserles, N{\o}rsett and Saff \cite{IsSa,IsNo,INS} were probably the first to point out the key role
of the real zeros of the polynomials $\sum_{j=0}^K \gamma_j(x)_j$ or $\sum_{j=0}^K \gamma_jx^j$ to prove the real rootedness of the polynomial
$$
\sum_{j=0}^K\gamma_j\rho_j L^\alpha _j(x),
$$
where $(\rho_j)_j$ are certain normalization sequences. They also proved similar results for many other sequences of orthogonal polynomials.
Iserles, N{\o}rsett and Saff results are related to the following problem: for which sequences of orthogonal polynomials $(p_n)_n$ with respect to a positive measure in the real line, the linear operator $T(x^n)=p_n$ preserves real-rootedness? (see \cite{piot,cha,CSW}). The definition of preserving real-rootedness is as follows.

\begin{Def}\label{prr}
Given a linear operator $T$ acting in the linear space of real polynomials such that $\deg T(p)=\deg p$, $T$ is said to preserve real-rootedness up to degree $K$ if for any polynomial $p$ having only real zeros and degree at most $K$ the polynomial $T(p)$ has also only real zeros. We say that $T$ preserves real-rootedness if $T$ preserves real-rootedness for all positive integer $K$.
\end{Def}

It is known that the operator $T(x^n)=p_n$ preserves real-rootedness when $p_n$ are either the Hermite or the Laguerre polynomials with parameter $\alpha=0$ (see \cite[p. 560]{IsSa} and \cite[Th. 1.2]{Durh} for the Hermite case and \cite{Fis} for the Laguerre case with $\alpha=0$). As a consequence of our results, we deduce the following.

\begin{Cor}\label{coj} The linear operator $T(x^n)=p_n$ preserves real-rootedness in the following cases:
\begin{enumerate}
\item If $p_n=L_n^\alpha$ and $\alpha>-1$ is not a non-negative integer then $T$ preserves real-rootedness up to degree $K$ if and only if $\alpha> K-2$.
\item If $p_n=L_n^\alpha$, with $\alpha$ a non-negative integer, then $T$ preserves real-rootedness.
\item If $p_n=\La_n^\alpha$ with $\alpha>-1$, then $T$ preserves real-rootedness.
\item If $p_n=\Ll_n^\alpha$ with $\alpha>-1$, then $T$ preserves real-rootedness.
\end{enumerate}
\end{Cor}
\medskip

Finally, in Section \ref{bez} we establish the asymptotic behaviour of the zeros of the polynomials $(q_n)_n$ for the four normalizations of the Laguerre polynomials considered in this paper.

\section{Preliminaries}
Along this paper, the interlacing property  is defined as follows.

\begin{Def}\label{dip} Given two finite sets $U$ and $V$ of real numbers ordered by size, we say that $U$ strictly interlaces $V$ if $\min U<\min V$ and between any two consecutive elements of any of the two sets there exists one element of the other.
If $U\cap V=W\not =\emptyset$, then we say that $U$ interlace $V$ if $U\setminus W$ strictly interlace $V\setminus W$.
\end{Def}
Observe that if $U$ interlaces $V$, then either $\card(U) = \card(V)$, and then $\max U < \max V$, or $\card(U) = 1 + \card(V)$, and then $\max U > \max V$. Observe also that the interlacing property is not symmetric, due to the condition $\min U < \min V$.

We will use the following version of Obreshkov theorem. Since we have not found that version in the literature, we include a proof (see, for instance, \cite[p. 9]{hou} for an standard version of Obreshkov theorem).

\begin{Theo}\label{obre}
Let $p$ and $q$ be real polynomials. Assume that $p$ has exactly $s$ real zeros, $q$ has exactly $s-1$ real zeros and all the zeros are simple.
If for all real $\lambda $ the real zeros of the polynomial $p(z)+\lambda q(z)$ are simple, then the real zeros of $p$ interlace the real zeros of $q$ (if in addition we assume that $p$ and $q$ have no zeros in common
then the interlacing is strict).

Assuming $\deg p=s$ and $\deg q=s-1$, then the converse is also true, i.e., if the zeros of $p$ interlace the zeros of $q$, then for all real $\lambda $ the polynomial $p(z)+\lambda  q(z)$ has only real and simple zeros.
\end{Theo}

\begin{proof}
By removing the common zeros of $p$ and $q$, we can assume that they have no common zeros, in which case, we will prove that the interlacing is strict.

Denote by $\xi_i$, $i=1,\cdots, s$, the real zeros of $p$ arranged in increasing order. We have to prove that $q$ has exactly one zero in each interval $(\xi_i,\xi_{i+1})$. For $x\in \RR$, define the polynomial
$$
h(z)=q(x)p(z)-p(x)q(z).
$$
Obviously $h(x)=0$, and since $p$ and $q$ have no common zeros, we deduce from the hypothesis that $x$ has to be a simple zero of $h$. Hence
$$
0\not =h'(x)=q(x)p'(x)-p(x)q'(x).
$$
This gives that $q(x)p'(x)-p(x)q'(x)$ has constant sign. If we evaluate at $\xi_i$, we conclude $q(\xi_i)p'(\xi_i)$ has constant sign. Since $p'(\xi_i)$ alternates sign with $i$, we deduce that $q$ has just one zero in each interval $(\xi_i,\xi_{i+1})$ (because $q$ has exactly $s-1$ real zeros).

If we assume that the zeros of $p$ interlace the zeros of $q$, and $s=\deg p=1+\deg q$, then $q(\xi_i)$ alternates sign with $i$.
We have then for the polynomial $h(z)=p(z)+\lambda q(z)$  that $h(\xi_i)=\lambda q(\xi_i)$, and so $h(\xi_i)$ alternates sign with $i$. Hence $h$ has a odd number of zeros in each interval $(\xi_i,\xi_{i+1})$. Since $\deg h=s$ we conclude that $h$ has $s-1$ simple zeros in $(\xi_1,\xi_s)$ and one more real zero less that $\xi_1$ or larger that $\xi_s$.

\end{proof}

The following elementary Lemmas will be useful (they are Lemmas 2.2 and 3.5 of \cite{Dur0}, respectively).

\begin{Lem}\label{lem2} Define from the numbers $A_j$, $j=0,\cdots , K$, $A_0,A_K\not =0$, the polynomial $P_A$ as
$$
P_A(x)=\sum_{j=0}^KA_j x^{K-j}.
$$
If  $\theta $ is a zero of $P_A$, we define the numbers $B_j$, $j=0,\cdots , K-1$, and the polynomial $P_{B}$ as
$$
P_{B}(x)=\frac{P_{A}(x)}{x-\theta}=\sum_{j=0}^{K-1}B_j x^{K-1-j}.
$$
Then, on the one hand, we have
\begin{equation}\label{id1l2}
A_j=\begin{cases} B_j-\theta B_{j-1},& j=1,\cdots, K-1,\\ B_0,& j=0,\\
-\theta B_{K-1},&j=K.\end{cases}
\end{equation}
And on the other hand
\begin{equation}\label{id2l2}
B_{j}=\sum_{i=0}^{j}\theta^iA_{j-i},\quad 1\le j\le K-1.
\end{equation}
\end{Lem}

\begin{Rem} \label{exp}
Using the notation of Lemma \ref{lem2}, given real numbers $B_j$, $0\le j\le K-1$, and $\theta$, we can produce the real numbers $A_j$, $0\le j\le K$, as in (\ref{id1l2}), so that we have for
$$
P_{A,\theta}=\sum_{j=0}^{K}A_jx^{K-j},\quad P_{B}(x)=\sum_{j=0}^{K-1}B_jx^{K-1-j}
$$
the identity
$$
P_{A,\theta}(x)=(x-\theta)P_{B},
$$
(we have included the real number $\theta$ in the notation $P_{A,\theta}$ to stress the dependence of this polynomial on $\theta$).

If we define
$$
q_n^{B}(x)=\sum_{j=0}^{K-1}B_jp_{n-j},\quad q_n^{A;\theta} (x)=\sum_{j=0}^{K}A_jp_{n-j}
$$
the identity (\ref{id1l2}) straightforwardly gives
\begin{equation}\label{idco}
q_n^{A;\theta}(x)=q_n^{B}(x)-\theta q_{n-1}^{B}(x),
\end{equation}
\end{Rem}

We then have.

\begin{Lem}\label{enze}
Assume that all the zeros of the polynomials $q_n^{B}$ are real and simple for $n\ge n_0$. Then the following conditions are equivalent.
\begin{enumerate}
\item The zeros of $q_{n+1}^{B}$ interlace the zeros of $q_n^{B}$ for $n\ge n_0$.
\item For all real number $\theta$ the polynomial $q_n^{A;\theta}$ has only real and simple zeros for $n\ge n_0+1$.
\end{enumerate}
Moreover, in that case the zeros of $q_n^{A;\theta}$ interlace the zeros of $q_n^B$.
\end{Lem}

The following Theorem due to Beardon and Driver will be also useful.

\begin{Theo}\cite[Theorem 3.1]{BDr}\label{BeDr} Let $(p_n)_n$ be orthogonal polynomials with respect to a positive measure. Fix $0<r<n$ and let $\zeta_i(n)$, $i=1,\dots, n$, the zeros of $p_n$ listed in increasing order. Let $P$ be a polynomial in the span of $p_r,\dots, p_n$. Then at least $r$ of the intervals $(\zeta_i,\zeta_{i+1})$ contain a zero of $P$.
\end{Theo}

In Section \ref{apen}, we need the concept of multiplier sequences.

\begin{Def}\label{mul} A sequence $T=(\lambda_n)_n$ of real numbers is called a multiplier sequence if, whenever the polynomial $p(x)=\sum_{j=0}^na_jx^j$ has only real zeros, the polynomial
$$
T(p(x))=\sum_{j=0}^n\lambda_j a_jx^j
$$
has only real zeros as well.
\end{Def}

Multipliers can be characterized using functions in the  Laguerre-Pólya class.

\begin{Def}\label{def1} An entire function $A$ is said to be in the Laguerre-Pólya class if it can be expressed in the form
\begin{equation}\label{pspr}
A(z)=cz^me^{-az^2+bz}\prod_{j=1}^\infty \left(1-\zeta_j z\right)e^{\zeta_j z},
\end{equation}
where $a\ge 0$, $m\in \NN$, $b, c,\zeta_j\in\RR$, $j\ge 1$, and $\sum_{j=1}^\infty \zeta_j^2<+\infty$. We say that an entire function $f$ is of type I (or first type) in the Laguerre-Pólya class, in short $f\in \lp I$, if
$f(z)$ or $f(-z)$ has a product representation of the form
$$
cz^me^{\alpha z}\prod_{j=1}^\infty \left(1+\zeta_jz\right),
$$
where $\alpha \ge 0$, $c\in \RR$, $m\in \NN$ and $\zeta_j> 0$, $\sum_{j=1}^\infty\zeta_j<\infty$.
\end{Def}

We then have (see \cite{PS}, \cite{CrCs}).

\begin{Theo}\label{PS2} A sequence $T=(\lambda_n)_n$ of real numbers is a multiplier if and only if the function
$$
f(z)=\sum_{n=0}^\infty \frac{\lambda_n}{n!}z^n
$$
is of type I in the Laguerre-Pólya class.
\end{Theo}

Examples of multipliers are the sequences $(1/(1+\alpha)_n)$, $\alpha>-1$.

\begin{Rem}\label{erdr} When necessary, we will use superscripts to stress the dependence of the sequences $q_n$ (\ref{qnmli}), (\ref{qnml5}), (\ref{qnmlii}) or (\ref{qnml3}) on the numbers $\gamma_j$, $1\le j\le K$, and the Laguerre parameter $\alpha$. For instance, we will sometimes  denote them $q_n^\gamma$ or $q_n^{\alpha;\gamma}$.
\end{Rem}

The following Lemma will be also useful (the proof is similar to the usual proof for the Hurwitz's Theorem (see \cite[p. 178]{ahl}) and it is omitted).

\begin{Lem}\label{mamon}
Let $f_n$, $g_n$, $f$ and $\eta$ be analytic functions in a region $\Omega$ of the complex plane. Assume that $f$ has $N$ non-real zeros in $\Omega$ and that
$$
\lim_nf_n(z)=\eta(z)f(z),\quad  \lim_ng_n(z)=f(z),
$$
uniformly in compact sets of $\Omega$.
Then there exists $n_*\in \NN$ such that for $n\ge n_*$ the function $f_n(z)+\lambda g_n(z)$ has at least $N$ non-real zeros in $\Omega$ for any real number $\lambda$.
\end{Lem}

(We stress that the positive integer $n_*$ guaranteed by the Lemma does not depend on the real number $\lambda$).

\section{Monic Laguerre polynomials}\label{ssec4}
In this Section, we prove  Theorem \ref{ldlc} by using the backward shift operator for the Laguerre polynomials
\begin{equation}\label{opl}
\Lambda_\alpha f(x)=xf'(x)+(\alpha-x)f(x),
\end{equation}
which it satisfies
\begin{equation}\label{opll}
\Lambda_\alpha (\hat L^\alpha_n)=-\hat L^{\alpha-1}_{n+1}.
\end{equation}

\begin{Def} A linear operator $T$ acting in the linear space of polynomials is a real zero increasing operator if for all polynomial $p$ the number of real zeros of $T(p)$ is greater than the number of real zeros of $p$.
\end{Def}

\begin{Lem} For $\alpha>0$, the operator $\Lambda_\alpha$ is a real zero increasing operator. Moreover, if all the zeros of $p$ are positive and simple then all the zeros of $\Lambda(p)$ are also positive and simple, and interlace the zeros of $p$.
\end{Lem}

\begin{proof}
Write $q=\Lambda_\alpha  p=xp'(x)+(\alpha-x)p(x)$.
Assume first that $p(0)\not =0$ and $p$ has $k$ real and simple zeros, say
$$
\zeta_1<\dots<\zeta_k.
$$
Assume in addition that all the real zeros are positive. Since $q(\zeta_i)=\zeta_ip'(\zeta_i)$, we deduce that $q$ has at least one zero in each interval $(\zeta_i,\zeta_{i+1})$, $i=1,\dots, k-1$; that is, at least $k-1$ positive zeros. Since the leading coefficient of $q$ and $p$ have different sign, we conclude that $q$ has to have a zero in the interval $(\zeta_k,+\infty)$, and then $q$ has at least $k$ real zeros. Using that the number of real zeros of a polynomial with real coefficients has to have the same parity as the degree, we deduce that $q$ has to have at least $k+1$ real zeros. Since $q(0)=\alpha p(0)$ implies that $q$ and $p$ has the same sign at $x=0$, we can conclude that $q$ has to have a zero in the interval $(0,\zeta_1)$; that is at least $k+1$ positive zeros. In particular, we have proved that if all the zeros of $p$ are positive and simple then all the zeros of $\Lambda(p)$ are also positive and simple, and interlace the zeros of $p$.

Assume next that $p$ has negative and positive zeros (and $p(0)\not =0$). More precisely $\zeta_l<0<\zeta_{l+1}$, for some $l$. Proceeding as before, we conclude that $q$ has at least one zero in each interval $(\zeta_i,\zeta_{i+1})$, $i=1,\dots, k-1$, $i\not =l$. That is, at least $k-2$ real zeros. If $p'(\zeta_l)<0$, we deduce that $p(0)<0$ and so $q(0)=\alpha p(0)<0$; since $q(\zeta_l)=\zeta_lp'(\zeta_l)>0$ and $q(\zeta_{l+1})=\zeta_{l+1}p'(\zeta_{l+1})>0$, we deduce that $q$ has at least one zero in each interval $(\zeta_l,0)$ and $(0,\zeta_{l+1})$. Hence $q$ has at least $k$ zeros, and as before, this implies that $q$ has to have al least $k+1$ zeros.

If $p$ has zeros of multiplicity bigger than $1$ or $p(0)=0$, we can use a continuity argument.

\end{proof}

Using the operator $\Lambda_\alpha$ we will find sufficient conditions on the real numbers $\gamma_j$, $1\le j\le K$, $\gamma_1=1$, $\gamma_K\not =0$, so that the polynomials
\begin{equation}\label{qnml}
q_n(x)=\sum_{j=0}^K \gamma_j\hat L^\alpha_{n-j}(x)
\end{equation}
have only real zeros for $n\ge K$.

In order to do that we consider the following polynomial
\begin{equation}\label{polpl}
Q(x)=\sum_{j=0}^K(-1)^j\gamma_{j}(x)_{K-j}.
\end{equation}

We next introduce some families of polynomials using the backward shift operator for the Laguerre polynomials (\ref{opl}). We will prove nice zero properties for the polynomials in these families, from where Theorem \ref{ldlc} will follow as corollary, because it turns out that the sequence of polynomials $(q_n)_n$ (\ref{qnml}) is a particular case of such families. This approach will also shows why the polynomial $Q$ is defined using the basis $((x)_n)_n$ of the linear space of polynomials.

Given a real number $r$ and two sequences of real numbers $\phi=(\phi_i)_{i\ge1}$ and $\psi=(\psi_i)_{i\ge1}$, we associate to $\phi, \psi$ the following sequence of polynomials: $\lge_0^{r;\phi,\psi}=1$ and for $n\ge 1$
\begin{equation}\label{dhl}
\lge_{n+1}^{r;\phi,\psi}(x)=\Lambda_{r}(\lge_{n}^{r;\phi,\psi}(x)) +(\phi_{n+1}+x\psi_{n+1})\lge_{n}^{r;\phi,\psi}(x).
\end{equation}

When $\psi_i\not =1$, for all $i\ge 1$,  then $\lge_{n}^{r;\phi,\psi}$ has degree $n$ and leading coefficient $\prod_{i=1}^n(\psi_i-1)$. From the definition, it follows easily that
\begin{equation}\label{eapm2}
\lge_{n}^{r+s;\phi,\psi}(x)=\lge_{n}^{r;s+\phi,\psi}(x),
\end{equation}
where the sequence $s+\phi$ is defined by $(s+\phi)_i=s+\phi_i$.

Under mild conditions on the parameters $r, \phi, \psi$, we can prove that the zeros of the polynomials $\lge_{n}^{r;\phi,\psi}$ behave nicely.

\begin{Theo}\label{pepl} Let $\phi$ and $\psi$ two sequences of real numbers. Given a positive integer $m$, assume that $r+\phi_i>0$ and $\psi_i<1$, $1\le i\le m$. Then for $n\le m$ the polynomial $\lge_{n}^{r;\phi,\psi}$ has only positive and simple zeros. Moreover, the zeros of $\lge_{n}^{r;\phi,\psi}$ strictly interlace the zeros of $\lge_{n-1}^{r;\phi,\psi}$. If, in addition, $\psi_{m+1}<1$ then the polynomial $\lge_{m+1}^{r;\phi,\psi}$ has only real and simple zeros (and at least $m$ of them are positive).
\end{Theo}

\begin{proof}
We proceed by induction on $n$. For $n=0,1$, the result is trivial.

Let $n\le m$ and assume $\lge_{n-1}^{r;\phi,\psi}$ has only positive and simple zeros. Write then
$$
0<\zeta_1<\dots<\zeta_{n-1},
$$
for the zeros of $\lge_{n-1}^{r;\phi,\psi}$.

The definitions (\ref{opl}) and (\ref{dhl}) give
$$
\lge_{n}^{r;\phi,\psi}(x)=x(\lge_{n-1}^{r;\phi,\psi})'(x) +[r+\phi_{n}+(-1+\psi_{n})x]\lge_{n-1}^{r;\phi,\psi}(x).
$$
And so
$$
\lge_{n}^{r;\phi,\psi}(\zeta_i)=\zeta_i(\lge_{n-1}^{r;\phi,\psi})'(\zeta_i).
$$
This shows that $\lge_{n}^{r;\phi,\psi}$ has an odd number of zeros in each interval $(\zeta_i,\zeta_{i+1})$, $i=1,\dots, n-2$.

Since the leading coefficient of $\lge_{n}^{r;\phi,\psi}$ is $\prod_{i=1}^{n}(\psi_i-1)$ and $\psi_i<1$, we deduce that the leading coefficients of $\lge_{n}^{r;\phi,\psi}$ and $\lge_{n-1}^{r;\phi,\psi}$ have different sign, and so the polynomial $\lge_{n}^{\phi,\psi}$ has also a number of odd zeros in the intervals $(-\infty,\zeta_1)$ and $(\zeta_{n-1},+\infty)$.
But $\lge_{n}^{r;\phi,\psi}(0)=(r+\phi_{n})\lge_{n-1}^{r;\phi,\psi}(0)$, and $r+\phi_{n}>0$, and hence we deduce that $\lge_{n}^{\phi,\psi}$ has a zero in the interval
$(0,\zeta_{1})$. This shows that $\lge_{n}^{\phi,\psi}$ has $n$ positive and simple zeros, and the interlacing property.

Proceeding in the same way, we can also prove that if $\psi_{m+1}<1$ then $\lge_{m+1}^{r;\phi,\psi}$ has only real and simple zeros.

\end{proof}

Given complex numbers $\phi_i$, $\psi_i$, $i\ge 1$, we can also define the polynomials $\lge _n^{r;\phi,\psi}$ as in (\ref{dhl}). The following result can be proved proceeding as in the proof of Theorem \ref{pepl}.

\begin{Cor}\label{cze} Let $\phi$ and $\psi$ two sequences of complex numbers.
Assume that $\lge _n^{r;\phi,\psi}$ has only positive and simple zeros for some $n$, and that $\phi_{n+1}$ and $\psi_{n+1}$ are real with $\psi_{n+1}<1$, then
$\lge _{n+1}^{r;\phi,\psi}$ has only real and simple zeros (and at least $n$ of them are positive). If, in addition, we also assume  that $r+\phi_{n+1}>0$, then
$\lge _{n+1}^{r;\phi,\psi}$ has only positive and simple zeros.
In both cases, the zeros of $\lge _{n+1}^{r;\phi,\psi}$  interlace the zeros of $\lge _{n}^{r;\phi,\psi}$.
\end{Cor}
\bigskip

When defining the polynomials $\lge_{n}^{r;\phi,\psi}$, we have avoided to call them \textit{generalized Laguerre} polynomials because
we do not recover the Laguerre polynomials when $\phi_i=\psi_i=0$. Indeed, when $\phi_i=\psi_i=0$ and $r=0$, we have
$$
\lge_{n}^{0;0,0}(x)=\Be_{n}(-x),
$$
where $\Be_n$ denotes the $n$-th Bell polynomial. These polynomials are defined by
\begin{equation}\label{comp1}
\Be_n(x)=\sum_{j=0}^n S(n,j)x^j,
\end{equation}
where $S(n,j)$, $0\le j\le n$, denote the Stirling numbers of the second kind defined in turn as
\begin{equation}\label{comp0}
x^n=\sum_{j=0}^nS(n,j)x(x-1)\cdots (x-j+1).
\end{equation}

The connection with the polynomials $\lge_{n}^{r;\phi,\psi}$ is because Bell polynomials satisfy the recurrence
$$
\Be_{n+1}(x)=x\left(1+\frac{d}{dx}\right)\Be_n(x).
$$
When $r\not =0$, we have
$$
\lge_{n}^{r;0,0}(x)=\Be_{n,r}(-x),
$$
where $\Be_{n,r}$ denotes the $n$-th $r$-Bell polynomial (see, for instance, \cite{mez,MeCo,MeRa}). We also have
$$
\lge_{n}^{r;0,s}(x)=\Be_{n,r}((s-1)x).
$$
In fact, the case $\psi_i=0$ has already been considered in \cite{durgb}, where the polynomials were called generalized Bell polynomials.

Laguerre polynomials appear when $r=\alpha+1$, $\phi_i=i-1$ and $\psi_i=0$ (it is an easy consequence of (\ref{opll})):
\begin{equation}\label{lqn}
\lge_{n}^{\alpha+1;\phi,0}(x)=(-1)^n\hat L^\alpha_{n}(x).
\end{equation}
This somehow unexpected connection between Laguerre and Bell polynomials (and taking into account (\ref{comp1}) and (\ref{comp0})) explains why the polynomial $Q$ (\ref{polpl}) is defined using the basis $((x)_n)_n$ of the linear space of polynomials.

In what follows we only consider the case when $\psi_i=0$ which, to simplify the notation, we will write $\lge_{n}^{r;\phi,0}=\lge_{n}^{r;\phi}$. In terms of the generalized Bell polynomials $(\lge_n^{\phi})_n$ introduced in \cite{durgb}, we have
$$
\lge_{n}^{r;\phi}(x)=\lge_{n}^{r+\phi}(-x).
$$
Notice that the polynomial $\lge_n^{r;\phi}(x)$ has degree $n$, leading coefficient equal to $(-1)^n$, and only depends on the numbers $\phi_1,\cdots, \phi_n$. Moreover, $\lge_n^{r;\phi}$ is a symmetric function of $\phi_1,\cdots, \phi_n$ (because the corresponding generalized Bell polynomial $\lge_{n}^{\phi}$ is so, see \cite[(1.6)]{durgb}).

\medskip

Write $\theta_i$, $i=1,\cdots, K$, for the zeros of the polynomial $Q$ (\ref{polpl})
and define the sequence $\phi$ as
\begin{equation}\label{sphl}
\phi_i=\begin{cases} -\theta_i,&1\le i\le K,\\
K-i,&K+1\le i.
\end{cases}
\end{equation}

We then have the following result.

\begin{Lem}\label{joe} Let $\alpha$ and $K$ be a real number and a positive integer, respectively. For numbers $\gamma_j$, $0\le j\le K$, $\gamma_0=1$, $\gamma_K\not =0$, consider the polynomial
$Q$ defined by (\ref{polpl}).
Then for all $n\ge K$, we have
\begin{equation}\label{idl}
q_n(x)=(-1)^n\lge_n^{\alpha+n-K+1;\phi}(x),
\end{equation}
where the polynomials $q_n$ and the sequence $\phi$ are defined in (\ref{qnml}) and (\ref{sphl}), respectively.
\end{Lem}

\begin{proof}
Define $A=\{\phi_i:1\le i\le K\}\setminus \{i-1: 1\le i\le K\}$ and write $N_A$ for the number of elements of $A$. We proceed by induction on $N_A$.

If $N_A=0$, we can assume
\begin{equation}\label{eapm3}
\phi_i=K-i, \quad 1\le i\le K
\end{equation}
(because the definition of $\lge_n^{r;\phi}$ only depends on the first $n$ elements of $\phi$, and the order does not matter since they are symmetric functions of $\phi_i$, $1\le i\le n$).
On the one hand, since $-\phi_i$, $1\le i\le K$, are the zeros of $Q$, we deduce that $Q(x)=x(x+1)\cdots (x+K-1)$. Hence from the definition of $Q$ (\ref{polpl}) we deduce $\gamma_0=1$ and $\gamma_j=0$, $1\le j\le K$, and so $q_n(x)=\hat L_n^\alpha(x)$.

On the other hand, using (\ref{eapm2}),we have
\begin{equation}\label{eapm1}
\lge_n^{\alpha+n-K+1;\phi}(x)=\lge_n^{\alpha+1;n-K+ \phi}(x).
\end{equation}
Writing $\tilde \phi_i=n-K+\phi_i$ and taking into account (\ref{sphl}) and (\ref{eapm3}), we have
$$
\tilde \phi_i=n-i,\quad i\ge 1.
$$
Since the first $n$ elements of $\tilde \phi$ are $\{0,1,\cdots, n-1\}$, we deduce using
(\ref{lqn}) and (\ref{eapm1}), that $\lge_n^{\alpha+n-K+1;\phi}(x)=(-1)^n\hat L_n^\alpha (x)$, for $n\ge K$. This proves the identity (\ref{idl}) when $N_A=0$.

Assume next that the identity (\ref{idl}) holds for $N_A=l-1$. We next prove it for $N_A=l$. Consider first the case $n=K$. Take $m$, $1\le m\le K$, such that $\phi_m\in A$ and consider the sequence $\phi^{\{m\}}$ defined by
$$
\phi^{\{ m\} }_i=\begin{cases} \phi_i, &1\le i\le m-1,\\
\phi_{i+1}, &m\le i,\end{cases}
$$
that is, $\phi^{\{ m\} }$ is que sequence obtained by removing the term $\phi_m$ from $\phi$.
Since the polynomial $\lge_K^{\alpha+1;\phi}$ only depends on $\phi_i$, $1\le i\le K$, we have from (\ref{dhl})
\begin{equation}\label{fi1}
\lge_K^{\alpha+1;\phi}(x)=\Lambda_{\alpha+1}\lge_{K-1}^{\alpha+1;\phi^{\{m\}}}(x)+\phi_m\lge_{K-1}^{\alpha+1;\phi^{\{m\}}}(x).
\end{equation}
Write
\begin{align}\label{fi2}
Q_{\phi_m}(x)&=\frac{Q(x)}{x-\phi_m}=\sum_{j=0}^{K-1}(-1)^j\tilde\gamma_j(x)_{K-1-j}.
\end{align}
Since $Q(x)=(x-\phi_m)Q_{\phi_m}(x)$ and
$$
(x-\phi_m)(x)_{K-j-1}=(x)_{K-j}-(K-j+\phi_m-1)(x)_{K-j-1},
$$
a simple computations gives
\begin{equation}\label{fi4}
\gamma_i=\tilde \gamma_i+(K-j-\phi_m)\tilde \gamma_{i-1},
\end{equation}
where we set $\tilde \gamma_K=\tilde \gamma_{-1}=0$.

The zeros of $Q_{\phi_m}$ are $-\phi_i$, $1\le i\le K$, $i\not =m$. We can then use the induction hypothesis to conclude that
$$
\lge_{K-1}^{\alpha+1;\phi^{\{m\}}}(x)=(-1)^{K-1}\sum_{j=0}^{K-1}\tilde\gamma_j\hat L_{K-1-j}^{\alpha}(x).
$$
Hence, from (\ref{fi1}) we have
$$
\lge_K^{\alpha+1;\phi}(x)=(-1)^{K-1}\left[\Lambda_{\alpha+1}\sum_{j=0}^{K-1}\tilde\gamma_j\hat L_{K-1-j}^{\alpha}(x)+\phi_m\sum_{j=0}^{K-1}\tilde\gamma_j\hat L_{K-1-j}^{\alpha}(x)\right].
$$
Using that $\hat L_{n}^{\alpha}=\hat L_{n}^{\alpha+1}+n\hat L_{n-1}^{\alpha+1}$, (\ref{opll}) and (\ref{fi4}), we get after easy computations
\begin{equation}\label{fi5}
\lge_K^{\alpha+1;\phi}(x)=(-1)^K\sum_{j=0}^{K}(\tilde\gamma_j+(K-j-\phi_m)\tilde\gamma_{j-1})\hat L_{K-j}^{\alpha}(x)=(-1)^K q_K(x).
\end{equation}
This is the identity (\ref{idl}) for $n=K$.

In order to prove the case $n=K+1$, we proceed as follows. On the one hand, we have already proved that
$$
\lge_K^{\alpha+2;\phi}(x)=(-1)^K\sum_{j=0}^{K}\gamma_j\hat L_{K-j}^{\alpha+1}(x).
$$
If we apply the operator $\Lambda_{\alpha+2}$ to that identity, we get using (\ref{opll})
$$
\Lambda_{\alpha+1} \lge_K^{\alpha+2;\phi}(x)=(-1)^{K+1}\sum_{j=0}^{K}\gamma_j\hat L_{K+1-j}^{\alpha}(x).
$$
On the other hand, taking into account that $\phi_{K+1}=-1$ (see (\ref{sphl})) and (\ref{eapm2}), the definition (\ref{dhl}) gives
\begin{align*}
\Lambda_{\alpha+1} \lge_K^{\alpha+2;\phi}(x)&=\Lambda_{\alpha+1} \lge_K^{\alpha+1;1+\phi}(x)\\&= \lge_{K+1}^{\alpha+1;1+\phi}(x)-(1+\phi_{K+1})\lge_K^{\alpha+1;1+\phi}(x)
\\&=\lge_{K+1}^{\alpha+2;\phi}(x).
\end{align*}
Hence
$$
\lge_{K+1}^{\alpha+2;\phi}(x)=(-1)^{K+1}\sum_{j=0}^{K}\gamma_j\hat L_{K+1-j}^{\alpha}(x),
$$
which it is the identity (\ref{idl}) for $n=K+1$.

The proof for the rest of the cases can be completed proceeding similarly.

\end{proof}

\begin{Rem}\label{cspm}
In the previous lemma, from the numbers $\gamma_j$, $0\le j\le K$, (no matter if they are real or not) we have constructed the polynomial $Q=Q_{\gamma}$ with zeros $\theta_i$, $1\le i\le K$. Conversely, from numbers $\theta_i$, we can recover the numbers $\gamma_j$ such that $Q_\gamma(x)=\prod_{i=1}^K (x-\theta_i)$ as follows. Given $K$  numbers $\theta_i$, $1\le i\le K$, define $\Phi^K_j$, $1\le j\le K$, by
\begin{equation}\label{lenx}
(x-\theta_1)\cdots (x-\theta_K)=\sum_{j=0}^K\Phi^K_jx^{K-j},
\end{equation}
and the  numbers $\gamma_j$, $0\le j\le K$, by
\begin{equation}\label{len2}
\gamma_{j}=(-1)^j\sum_{i=0}^j(-1)^iS(K-j+i,K-j)\Phi^K_{j-i},
\end{equation}
where $S(n,i)$ denote the Stirling numbers of the second kind (see (\ref{comp0})).
Then
\begin{equation}\label{lenxy}
\prod_{i=1}^K (x-\theta_i)=\sum_{j=0}^K(-1)^j\gamma_j(x)_{K-j}.
\end{equation}

Indeed, using
$$
x^{K-j}=\sum_{i=0}^{K-j}S(K-j,i)x(x-1)\cdots (x-i+1)
$$
(see (\ref{comp0})) in (\ref{lenx}), we can deduce (\ref{lenxy}) after straightforward computations.

\end{Rem}

We are now ready to prove Theorem \ref{ldlc}.

\begin{proof}[Proof of Theorem \ref{ldlc}]
We define
$$
X_1=\{\theta\in\RR:\mbox{$Q(\theta)=0$ and $\theta\ge \alpha+1$}\}, \quad X_2=\{\theta:\mbox{$Q(\theta)=0$ and $\theta$ is non-real}\}.
$$
Write $r_1$ and $r_2$ for the number of elements of $X_1$ and $X_2$, respectively, and arrange the zeros $\theta_i$, $1\le i\le K$, of $Q$ such that $\theta _i\in X_1$ for $1\le i\le r_1$, $\theta _i\in X_2$ for $r_1+1\le i\le r_1+r_2$, and then $\theta_i<\alpha+1$ if $r_1+r_2+1 \le i\le K $.

We proceed in four steps.

\medskip

\noindent
\textit{Step 1.} If the polynomial $Q$ has only real zeros then all the zeros of the polynomial $q_n$ (\ref{qnml}) are real and simple for $n\ge n_1$, and positive and simple for $n\ge n_0$, where
$n_l= \max\{K,\lfloor\theta_i-\alpha+K\rfloor, 1\le i\le r_1-l\}$.
\begin{proof}[Proof of Step 1]
If we use (\ref{sphl}) and (\ref{idl}), we have
$$
\alpha+n-K+1+\phi_i=\begin{cases} \alpha+n-K+1-\theta_i,& 1\le i\le K,\\ \alpha+1+n-i,& K+1\le i\le n. \end{cases}
$$
Since $\theta_i$ are the zeros of $Q$, if we assume that $n\ge n_0$, we have that $\theta_i-\alpha-1+K\le n_0$, $1\le i\le r_1$, and
$\theta_i-\alpha-1+K<K\le n_0$, $r_1+1\le i\le K$. Hence, using that $\alpha+1>0$, we conclude that
$$
\alpha+n-K+1+\phi_i>0,\quad 1\le i\le n.
$$
Theorem \ref{pepl} then says that all the zeros of $q_n$ are positive and simple. The proof is similar if we assume $n\ge n_1$.
\end{proof}

\medskip

\noindent
\textit{Step 2.} If $Q$ has non-real zeros, then there exists a positive integer $n_0$, depending only on the zeros $\theta_j$, $1\le j\le r_1+r_2$, of $Q$,
such that all the zeros of the polynomial $q_n$  are positive and simple for $n\ge n_0$. Moreover, the positive integer $n_0$ depends uniformly on the zeros $\theta_j$, $1\le j\le r_1+r_2$, of $Q$.
\begin{proof}[Proof of Step 2]
Write $H=r_1+r_2$, and then for $1\le i\le H$, $\theta _i$ is either non-real or $\theta_i\ge \alpha+1$. For a positive integer $I$, $H\le I\le K$, consider the real numbers $\gamma_j^{(I)}$, $0\le j\le I$, defined by
\begin{equation}\label{teta}
\prod_{j=1}^I(x-\theta_j)=\sum_{j=0}^I(-1)^j\gamma_j^{(I)} (x)_{I-j}.
\end{equation}
Define then the sequence $\phi^{(I)}$, $H\le I\le K$, as follows
\begin{equation}\label{sphlo}
\phi_i^{(I)}=\begin{cases} -\theta_i,&1\le i\le I,\\
I-i,&I+1\le i.
\end{cases}
\end{equation}
The identity (\ref{idl}) gives for the sequence $\phi^{(I)}$ that
\begin{equation}\label{mvm}
q_{n}^{\gamma ^{(I)}}=(-1)^{n}\lge_{n}^{\alpha+n-I+1;\phi^{(I)}},\quad n\ge I.
\end{equation}
Corollary \ref{1lagi} implies the existence of a positive integer $n_1\ge K$ such that $q_{s}^{\gamma ^{(H)}}$ has only positive and simple zeros for $s\ge n_1$. This number $n_1$ only depends on the zeros $\theta_i$, $1\le i\le H$, of $Q$ (see the definition of $\gamma_j^{(I)}$, $0\le j\le I$, in (\ref{teta})). Moreover, this dependence is uniform of the zeros $\theta_i$, $1\le i\le H$.
This gives that the polynomial
$$
q_{s}^{\gamma ^{(H)}}=(-1)^{s}\lge_{s}^{\alpha+s-H+1;\phi^{(H)}},\quad s\ge n_1,
$$
has only positive and simple zeros.

For $s\ge n_1$, consider next the sequence $\phi^{s;(H+1)}$ defined by
$$
\phi_i^{s;(H+1)}=\begin{cases} -\theta_i,&1\le i\le H,\\
H-i,&H+1\le i\le s,\\  -\theta_{H+1},& i=s+1.
\end{cases}
$$
On the one hand, since $\theta_{H+1}<\alpha+1$, we have ($s\ge n_1\ge K\ge H$)
$$
\alpha+s-H+1+\phi_{s+1}^{s;(H+1)}=\alpha+s-H+1-\theta_{H+1}>s-H+1>0.
$$
Using the first part of Corollary \ref{cze}, we deduce that for $s\ge n_1$
the polynomial
$$
\lge_{s+1}^{\alpha+s-H+1;\phi^{s;(H+1)}}
$$
has only positive and real zeros.

On the other hand, since $\{\phi^{(H+1)}_i: 1\le i\le s+1\}=\{\phi^{s;(H+1)}_i: 1\le i\le s+1\}$, we have
$$
\lge_{s+1}^{\alpha+s-H+1;\phi^{s;(H+1)}}=\lge_{s+1}^{\alpha+s-H+1;\phi^{(H+1)}}
$$
(because the polynomial $\lge_n^{r;\phi}$ is a symmetric function of $\phi_i$, $1\le i\le n$). Using (\ref{mvm}), we deduce that
the polynomial
$$
q_{s+1}^{\gamma ^{(H+1)}}=(-1)^{s+1}\lge_{s+1}^{\alpha+s+1-(H+1)+1;\phi^{(H+1)}}
$$
has only positive and real zeros for $s\ge n_1$.

Proceeding similarly and taking into account that $\phi^{(K)}=\phi$ (\ref{sphlo}) and $\gamma^{(K)}=\gamma$, we deduce that
$$
q_{s+K-H}^{\gamma}=(-1)^{s+K-H}\lge_{s+K-H}^{\alpha+s+K-H-K+1;\phi}
$$
has only positive and simple zeros for $s\ge n_1$. That is
$$
q_{n}^{\gamma}=(-1)^{n}\lge_{n}^{\alpha+n-K+1;\phi}
$$
has only positive and simple zeros for $n\ge n_0=n_1+K-H$.

The proof also shows that $n_0$ depends uniformly on the zeros $\theta_j$, $1\le j\le r_1+r_2=H$, of $Q$.
\end{proof}

\medskip

So, we have proved parts (1) and (2) of Theorem \ref{ldlc}.

Even when the zeros of $q_{n+1}$ and $q_n$ are all real, they do not need to interlace. Here it is a counter example: $K=2$ and
$$
q_n(x)=\hat L_n^0(x)+\frac{14}5\hat L_{n-1}^0(x)+\frac{81}{100}\hat L_{n-2}^0(x),\quad n\ge 2.
$$
Hence, $Q(x)=(x+9/10)^2$ (see (\ref{polpl})). Step 1 gives then that the polynomial $q_n$
has only real zeros for $n\ge 2$. But the zeros of $q_{n+1}$ do not interlace the zeros of $q_n$ for $n=2,3$.

However, we can guarantee the interlacing property for $n$ big enough. In order to prove it we will use the Lemma \ref{enze}.
So, we write the polynomial $Q$ (\ref{polpl}), as $Q(x)=\sum_{j=0}^KB_jx^{K-j}$, and hence $B_0=1$.
Take next a real number $\upsilon$ and define the polynomials $P_B(x)=Q(x)$,
$$
P_A(x)=(x-\upsilon)Q(x)=\sum_{j=0}^{K+1}A_jx^{K+1-j}.
$$
Setting $q_n^A$ and $q_n^B$ as is (\ref{idco}), if we prove that there exists $n_0$ (which does not depend on $\upsilon$) such that $q_n^{A}$ has only real zeros for any $n\ge n_0$, then the interlacing property follows from Lemma \ref{enze}.

We proceed in two more steps.

\medskip

\noindent
\textit{Step 3.} There exist $M>0$ and a positive integer $n_1$ (which do not depend on $\upsilon$) such that, for $|\upsilon|>M$, all the zeros of $q_n^{A}$ are real and simple for $n\ge n_1$.
\begin{proof}[Proof of Step 3]
Indeed, define $Q_A(x)=\sum_{j=0}^{K+1}(-1)^jA_j(x)_{K+1-j}$. Since
$$
\lim_{\upsilon \to \infty} \frac{P_{A}(x)}{\upsilon}=-Q(x),
$$
we deduce that
$$
\lim_{\upsilon \to \infty} \frac{Q_{A}(x)}{\upsilon}=\sum_{j=0}^{K}(-1)^jB_j(x)_{K-j}=R(x),
$$
and $R$ is a monic polynomial of degree $K$  which does not depend on $\upsilon$.

This shows that there exists a zero $\theta^A(\upsilon )$ of $Q_A$ which goes to $\infty$ when $\upsilon\to \infty$, and the other zeros of $Q_A$ goes to the zeros of the polynomial $R$. Since the coefficients of $Q_A$ are real, $\theta^A(\upsilon )$ has to be real for $\upsilon$ big enough. We then write $\theta_i^A$, $1\le i\le K+1$, for the zeros of $Q_A$ where $\theta_{K+1}^A=\theta^A(\nu)$.
Define next the sequences $\phi^A$, $\tilde\phi^A$ as
$$
\phi_i^A=\begin{cases} -\theta_i^A,&1\le i\le K+1,\\
K+1-i,&K+2\le i.
\end{cases},\quad
\tilde \phi_i^A=\begin{cases} -\theta_i^A,&1\le i\le K,\\
K-i,&K+1\le i.
\end{cases}
$$
Define the real numbers $\gamma_j^A$, $1\le i\le K+1$, from $\theta_i^A$, $1\le i\le K+1$, as explained in the Remark \ref{cspm}, and, in the same way, $\tilde \gamma_j^A$, $1\le i\le K$, from $\theta_i^A$, $1\le i\le K$.
Consider finally the polynomials
$$
q_n^{\gamma^A}(x)=\sum_{j=0}^{K+1} \gamma_j^A\hat L_{n-j}^\alpha(x), \quad q_n^{\tilde \gamma^A}(x)=\sum_{j=0}^{K} \tilde \gamma_j^A\hat L_{n-j}^\alpha(x).
$$
In particular, by the construction, on the one hand, we have $\gamma_j^A=A_j$, $1\le i\le K+1$, and so $q_n^A=q_n^{\gamma^A}$. And on the other hand
$$
q_n^{\tilde \gamma^A}=(-1)^n\lge_n^{\alpha+n-K+1;\tilde\phi ^A},\quad q_{n+1}^{\gamma^A}=(-1)^{n+1}\lge_{n+1}^{\alpha+n+1-(K+1)+1;\phi ^A}.
$$
Since the zeros $\theta_i^A$, $1\le i\le K$, of $Q_A$ converge to the zeros of $R$ as $\upsilon\to \infty$, using the Step 2 above, we deduce the existence of $M>0$ and a positive integer $n_1$ (which do not depend on $\upsilon$) such that for $n\ge n_1$ and $|\upsilon|>M$, $\theta^A(\upsilon )$ is a real number and the polynomial
$$
q_n^{\tilde \gamma^A}=(-1)^n\lge_n^{\alpha+n-K+1;\tilde\phi ^A}
$$
has positive and simple zeros.

Since $\{\phi^{A}_i: 1\le i\le n+1\}=\{\theta^A(\upsilon )\}\cup \{\tilde\phi^{A}_i: 1\le i\le n\}$, and $\theta^A(\upsilon )$ is a real number, using  Corollary \ref{cze}, we deduce that
$$
\lge_{n+1}^{\alpha+n-K+1;\phi ^A}
$$
has real and simple zeros for $n\ge n_1$. The Step 3 follows taking into account that
$$
q_{n+1}^{A}=q_{n+1}^{\gamma^A}=(-1)^{n+1}\lge_{n+1}^{\alpha+n+1-(K+1)+1;\phi ^A}.
$$
\end{proof}

\medskip

\noindent
\textit{Step 4.} There exists a positive integer $n_2$ (which does not depend on $M$) such that $q_n^{A}$ has only positive and real zeros for $n\ge n_2$ and $|\upsilon|\le M$.
\begin{proof}[Proof of Step 4]
Since
$$
q_{n}^A(x)=\sum_{j=0}^{K+1} A_j \hat L_{n-j}^\alpha(x),
$$
Step 4 is an easy consequence of Corollary \ref{1lagi} and (\ref{id1l2}).
\end{proof}

Since neither $n_1$ in Step 3 nor $n_2$ in Step 4 depends on $\upsilon$, we have that
$n_0=\max\{n_1+1,n_2\}$ does not depend on $\theta$, and $q_n^{A}$ has only real and simple zeros for $n\ge n_0$.
\end{proof}

\bigskip

Remark \ref{cspm} allows us to generate  numbers $\gamma_j$, $j=0,\dots, K$, for which the polynomials $q_n$ (\ref{qnml}) has only positive zeros for all $n\ge K$, using the following version of Part (1) in
Theorem \ref{ldlc}.

\begin{Cor}\label{eda} Let $\alpha$ and $K$ be a real number $\alpha>-1$ and a positive integer, respectively.
Given $K$  real numbers $\theta_i$, $1\le i\le K$, satisfying $\theta_i<\alpha+1$, define $\Phi^K_j$
and $\gamma_{j}$, $1\le j\le K$, as in (\ref{lenx}) and (\ref{len2}), respectively.
Then all the zeros of the polynomial $q_n$ (\ref{qnml}) are positive and simple for $n\ge K$. For $n$ big enough the zeros of $q_{n+1}$ interlace the zeros of $q_n$.
\end{Cor}

\section{Laguerre polynomials taking the value one at zero}\label{secl33}
In this Section we normalize the Laguerre polynomials so that they take the value one at zero
\begin{equation}\label{lag2n}
\La_n^\alpha(x)=\frac{n!}{(1+\alpha)_n}L_n^\alpha(x).
\end{equation}
It is easy to check that
\begin{equation}\label{val0}
(\La_n^\alpha)^{j}(0)=\frac{(-1)^jn(n-1)\cdots (n-j+1)}{(\alpha+1)_j}.
\end{equation}

For real numbers $\gamma_j$, $j=0,\cdots, K$, $\gamma_0=1$ and $\gamma_K\not=0$, we define
\begin{equation}\label{qnu}
q_n^{\alpha}(x)=\sum_{j=0}^K\gamma_j\La_{n-j}^\alpha(x),\quad n\ge K,
\end{equation}
and
\begin{equation}\label{dPl2}
P(x)=\sum_{j=0}^K\gamma_jx^{K-j}.
\end{equation}

The following Lemma is an easy consequence of (\ref{val0})

\begin{Lem}\label{val00}
The following conditions are equivalent.
\begin{enumerate}
\item The polynomial $q_{n_0}$ has a zero at $x=0$ of multiplicity $s$ for some $n_0\ge K$.
\item The polynomial $P$ has a zero at $x=1$ of multiplicity $s$.
\item The polynomial $q_{n}$ has a zero at $x=0$ of multiplicity $s$ for all $n\ge K$.
\end{enumerate}
\end{Lem}

\medskip

Consider next the following first order differential operator ($\alpha\not =0$)
\begin{equation}\label{uu2}
\Upsilon_\alpha p(x)=p(x)+\frac{x}{\alpha} p'(x).
\end{equation}
It is easy to check that $\Upsilon_\alpha(\La_n^\alpha)=\La_n^{\alpha-1}$, $n\ge 0$. And hence
\begin{equation}\label{qgs}
\Upsilon_\alpha (q_n^{\alpha})=q_n^{\alpha-1}.
\end{equation}

We need to consider the following function $\vartheta$ which measures the total multiplicity of multiple real zeros of a real polynomial $p$. More precisely, let $\xi_j$ be, $1\le j\le u$, the zeros of $p$ with multiplicity $m_j>1$, $1\le j\le u$, then we define $\vartheta(p)=\sum_{j=1}^um_j$.

In the next Lemma, we prove that the operator $\Upsilon_\alpha$, $\alpha>0$, has nice properties preserving the real zeros of polynomials.

\begin{Lem}\label{nulu}
Let $p$ be a real polynomial and denotes by $Z(p)$ the set of real zeros of $p$. Assume $\alpha >0$. Then
\begin{enumerate}
\item $\Upsilon_\alpha$ is a complex zero decreasing operator (i.e., the number of non-real zeros of $\Upsilon_\alpha p$ is less than or equal to the number of non-real zeros of $p$).
\item If $c<d$ are two consecutive elements of $Z(p)\cup \{0\}$ then $\Upsilon_\alpha p$ has some zero in $(c,d)$.
\item If $p$ has exactly $N^+$ positive and $N^-$ negative zeros, then $\Upsilon_\alpha p$ has at least $N^+$ positive and $N^-$ negative zeros.
\item If $p$ and $\Upsilon_\alpha p$ has the same number of real zeros and $p(0)\not=0$, then $\vartheta(\Upsilon_\alpha p)<\vartheta(p)$.
\end{enumerate}
\end{Lem}

\begin{proof}
Since $\Upsilon _\alpha (x^j)=(\alpha+j)x^j$ using \cite[Part 3 of Theorem 2.4]{CrCs}, we deduce that $\Upsilon_\alpha$ is a complex zero decreasing operator.

The proof of part (2) is as follows. Assume first that $c,d\not =0$, that is, $c$ and $d$ are zeros of $p$. Hence we can write
$p(x)=(x-c)^n(d-x)^mr(x)$, where the polynomial $r$ has constant sign in $[c,d]$. This gives $q(x)=\alpha \Upsilon _\alpha p(x)=(x-c)^{n-1}(d-x)^{m-1}s(x)$, where
$$
s(x)=x[(n(d-x)-m(x-c))r(x)+(x-c)(d-x)r'(x)]+\alpha(x-c)(d-x)r(x).
$$
And hence
\begin{align*}
s(c)&=cn(d-c)r(c),\\
s(d)&=-dm(d-c)r(d).
\end{align*}
Since $c,d$ are consecutive elements of $Z(p)\cup \{0\}$ and $c,d\not=0$, we deduce that $c$ and $d$ has equal sign, and then $s(c)s(d)<0$. So $s$ has some zero in $(c,d)$ and so $ \Upsilon _\alpha p$ as well.

If $d=0$ and it is also a zero of $p$. We write as before $p(x)=(x-c)^n(-x)^mr(x)$, where the polynomial $r$ has constant sign in $[c,0]$. This gives $q(x)=\alpha \Upsilon _\alpha p(x)=(x-c)^{n-1}(-x)^{m}s(x)$, where
$$
s(x)=(nx+m(x-c))r(x)+x(x-c)r'(x)+\alpha(x-c)r(x).
$$
And hence
\begin{align*}
s(c)&=cnr(c),\\
s(0)&=-c(m+\alpha)r(0).
\end{align*}
Since $c<0=d$, then $s(c)s(0)<0$. So $s$ has some zero in $(c,d)$ and so $ \Upsilon _\alpha p$ as well.

If $d=0$ and $p(0)\not =0$, we can write $q(x)=(x-c)^nr(x)$, where the polynomial $r$ has constant sign in $[c,d]$. This gives $\alpha q(x)=\Upsilon _\alpha p(x)=(x-c)^{n-1}s(x)$, where
$$
s(x)=x(nr(x)+(x-c)r'(x))+\alpha(x-c)r(x).
$$
And hence
\begin{align*}
s(c)&=cnr(c),\\
s(0)&=-\alpha cr(0).
\end{align*}
Since $\alpha>0$, we deduce that $s(c)s(0)<0$. So $s$ has some zero in $(c,0)$ and so $ \Upsilon _\alpha p$ as well.

The case $c=0$ can be proved similarly.

Part (3) is an easy consequence of part (2).

We finally prove the Part (4). Since $p(0)\not =0$, it is easy to see that if $\xi$ is a zero of $p$ of multiplicity $m>1$ then $\xi$ is a zero of $\Upsilon _\alpha p$ of multiplicity $m-1$.
This gives that $\Upsilon _\alpha p$ has at least $\vartheta (p)-u$ zeros that are also zeros of $p$ (counting with multiplicity), where $u$ is the number of multiple zeros of $p$.
Denote by $l$ the number of real zeros of $p$ (counting again with multiplicity). Then the cardinal of the set $Z(p)\cup\{0\}$ is $l-\vartheta (p)+u+1$. Using that $\Upsilon _\alpha p$ has also $l$ real zeros, we deduce from Part (2) that $\Upsilon _\alpha p$ has exactly one zero between two consecutive elements of $Z(p)\cup\{0\}$, and so these zeros have to be simple. This gives
$\vartheta (\Upsilon_\alpha p)\le \vartheta (p)-u<\vartheta (p)$.
\end{proof}

We need to consider the conformal mapping $\varphi$  of $\CC\setminus [0,4]$ onto $\{w:|w|>1\}$ defined by
\begin{equation}\label{defp}
\varphi(z)=\frac{1}2(z-2+\sqrt{z^2-4z}), \quad z\in \CC\setminus [0,4].
\end{equation}
with $\sqrt{z^2-4z}>0$ when $z>4$.

\begin{Rem}\label{remw} Given $w\in \CC$ with $|w|>1$, it is easy to prove that there exists an unique $z\in \CC\setminus [0,4]$ such that $w=-\varphi(z)$. Moreover,
$$
z=2-w-\frac1w,
$$
and $z$ is non-real if and only if $w$ is non-real.
\end{Rem}

We then have the following asymptotic for the Laguerre polynomials (see \cite{JeVa}, p. 63): for $z\in \CC\setminus [0,4]$,
\begin{align}\label{eapm5}
\lim_{n\to \infty}\sqrt{2\pi n}(-1)^n\frac{L_n^\alpha(nz)}{(\varphi(z))^n}e^{\frac{-nz}{1+\varphi(z)}}
=\frac{\varphi(z)^{1/2}(1+\varphi(z))^\alpha}{z^\alpha(z^2-4z)^{1/4}}.
\end{align}
Using that $L_n^{\alpha-1}(x)=L_n^\alpha(x)-L_{n-1}^\alpha(x)$, and the asymptotic (\ref{eapm5}), we deduce for $s\in \NN\setminus\{0\}$ and $z\in \CC\setminus [0,4]$
\begin{equation}\label{eapm6}
\lim_{n\to \infty}\frac{L_{n-s}^\alpha(nz)}{L_{n}^\alpha(nz)}=\frac{(-1)^s}{(\varphi(z))^{s}}.
\end{equation}

Since
\begin{equation}\label{gan2}
\frac{n!}{(1+\alpha)_n}=\Gamma(1+\alpha)\frac{\Gamma(1+n)}{\Gamma(1+\alpha+n)}\sim \Gamma(1+\alpha)n^{-\alpha},
\end{equation}
when $n\to \infty$ (see \cite[vol. I (4), p. 47]{EMOT}), the asymptotic (\ref{eapm5}) gives: for $z\in \CC\setminus [0,4]$,
\begin{align}\label{eapu}
\lim_{n\to \infty}\sqrt{2\pi}(-1)^n n^{\alpha+1/2}\frac{\La_n^\alpha(nz)}{(\varphi(z))^n}e^{\frac{-nz}{1+\varphi(z)}}
&=\frac{\Gamma(1+\alpha)\varphi(z)^{1/2}(1+\varphi(z))^\alpha}{ z^\alpha(z^2-4z)^{1/4}},\\\label{eapuu}
\lim_{n\to \infty}\frac{\La_{n-s}^\alpha(nz)}{\La_{n}^\alpha(nz)}&=\frac{(-1)^s}{(\varphi(z))^{s}}.
\end{align}

\begin{Lem}\label{aux20} For $z\in \CC\setminus [0,4]$, we have the following asymptotic for the polynomials $q_n$ (\ref{qnu})
\begin{equation}\label{asqnu}
\lim_{n\to \infty}
\frac{(-1)^n\sqrt{2\pi }n^{\alpha+1/2} q_n(nz)}{(\varphi(z)e^{\frac{z}{1+\varphi(z)}})^n}
=\frac{(-1)^K\Gamma(1+\alpha)(1+\varphi(z))^\alpha}{\varphi(z)^{-1/2+K} z^\alpha(z^2-4z)^{1/4}}
P(-\varphi(z)).
\end{equation}
\end{Lem}

\begin{proof}
It is just a matter of calculation using (\ref{eapu}) and (\ref{eapuu}):
\begin{align*}
\lim_{n\to \infty}
\frac{(-1)^n\sqrt{2\pi }n^{\alpha+1/2}q_n(nz)}{(\varphi(z)e^{\frac{z}{1+\varphi(z)}})^n}&=\lim_{n\to \infty}
\frac{(-1)^n\sqrt{2\pi }n^{\alpha+1/2}\La_n^\alpha(nz)}{(\varphi(z)e^{\frac{z}{1+\varphi(z)}})^n}
\sum_{j=0}^K\gamma_j\frac{\La_{n-j}^\alpha(nz)}{\La_n^\alpha(nz)}\\\nonumber
&=\frac{\Gamma(1+\alpha)\varphi(z)^{1/2}(1+\varphi(z))^\alpha}{z^\alpha(z^2-4z)^{1/4}}
\sum_{j=0}^K\gamma_j\frac{(-1)^j}{(\varphi(z))^{j}}\\\nonumber
&=\frac{(-1)^K\Gamma(1+\alpha)(1+\varphi(z))^\alpha}{\varphi(z)^{-1/2+K}z^\alpha(z^2-4z)^{1/4}}
\sum_{j=0}^K\gamma_j(-\varphi(z))^{K-j}\\\label{asqn}
&=\frac{(-1)^K\Gamma(1+\alpha)(1+\varphi(z))^\alpha}{\varphi(z)^{-1/2+K}z^\alpha(z^2-4z)^{1/4}}
P(-\varphi(z)).
\end{align*}
\end{proof}

We are now ready to prove Theorem \ref{hithis}.

\begin{proof}[Proof of Theorem \ref{hithis}]

We proceed in six steps.

\medskip

\noindent
\textit{Step 1.} We can assume that $P(1)\not=0$.

\begin{proof}[Proof of Step 1]
Otherwise if $P$ has a zero at $x=1$ of multiplicity $s$, Lemma \ref{val00} says that then $q_n$ has a zero of multiplicity $s$ at $x=0$. Hence, if we take
$$
\tilde P(x)=P(x)/(x-1)^s=\sum_{j=0}^{K-s}\tilde \gamma_jx^{K-s-j},
$$
and
$$
\tilde q_n(x)=\sum_{j=0}^{K-s}\tilde \gamma_j\La_{n-j}^\alpha,
$$
we have $q_n(x)=x^s\tilde q_n(x)$, $\tilde P(1)\not =0$, and the properties for the zeros of $q_n$ would follow from the properties of the zeros of $\tilde q_n$.
\end{proof}

Hence, we assume $P(1)\not =0$.

\medskip

\noindent
\textit{Step 2.} For each $n\ge K$ there exists $\alpha_n>0$ such that
$q_n^\alpha(x)$ (see (\ref{qnu})) has to have exactly $N^{\nr}$ non-real zeros and exactly $N^1$ negative zeros
for $\alpha>\alpha_n$.

\begin{proof}[Proof of Step 2]
Using the definition of the Laguerre polynomials, we have
$$
\La_n^\alpha(\alpha x)=\frac{n!}{(1+\alpha)_n}\sum_{l=0}^n\binom{n+\alpha}{n-l}\frac{\alpha^l(-x)^l}{l!}.
$$
This shows that
$$
\lim_{\alpha\to +\infty}\La_n^\alpha(\alpha x)=(1-x)^n.
$$
We start considering the case $n=K$. Hence
$$
\lim_{\alpha\to +\infty}q_K^\alpha(\alpha x)=P(1-x).
$$
If all the zeros of $P$ are simple, this shows that for $\alpha$ big enough $q_K^\alpha(x)$ has to have exactly $N^{\nr}$ non-real zeros and exactly $N^1$ negative zeros.
Since $P(1)\not=0$, as a consequence of Lemma \ref{val00}, using an argument of continuity, we deduce that, even if some of the zeros of $P$ are not simple, for $\alpha$ big enough $q_K^\alpha(x)$ has to have exactly $N^{\nr}$ non-real zeros and exactly $N^1$ negative zeros.

For $n\ge K$, write
$$
\tilde \gamma_j=\begin{cases}\gamma_j,&0\le j\le K,\\ 0, &K+1\le j.\end{cases}.
$$
Hence, if we define
$$
P(x;n)=P(x)x^{n-K}=\sum_{j=0}^n\tilde\gamma_jx^{n-j},
$$
we have
$$
q_n^\alpha(x)=\sum_{j=0}^n\tilde \gamma_jL_{n-j}^\alpha(x).
$$
Since $P(1;n)=P(1)\not =0$, we deduce that there exists $\alpha_n>-1$ such that for $\alpha>\alpha_n$ the polynomial $q_n^\alpha(x)$ has exactly $N^{\nr}$ non-real zeros and exactly $N^1$ negative zeros.
\end{proof}

\medskip

\noindent
\textit{Step 3.} For  $\alpha>-1$ and $n\ge K$, the polynomial $q_n^\alpha$
has to have at most $N^{\nr}$ non-real zeros, at least $N^1$ negative zeros and at least $n-N^{\nr}-N^1$ positive zeros.

\begin{proof}[Proof of Step 3]
For $\alpha>-1$ take $u\ge n$ such that $\alpha+u>\alpha _n$. This means that the polynomial $q_n^{\alpha+u}(x)$ has $N^{\nr}$ non-real zeros and $n-N^{\nr}$ real zeros of which $N^1$ are negative and $n-N^{\nr}-N^1$ are positive.
Since $\alpha+u-j>0$, $0\le j\le u-1$, we can use (\ref{qgs}) and Parts (1) and (3) of Lemma \ref{nulu} to conclude that the polynomial
\begin{equation}\label{cadu}
\Upsilon_{\alpha+1}\circ\cdots\circ\Upsilon_{\alpha+u}(q_n^{\alpha+u})=q_n^\alpha
\end{equation}
has to have at most $N^{\nr}$ non-real zeros, at least $N^1$ negative zeros and at least $n-N^{\nr}-N^1$ positive zeros.
\end{proof}

\medskip

\noindent
\textit{Step 4.} If all the zeros of $P$ are real, then for $n\ge K$, the polynomial $q_n^\alpha$ has only real and simple zeros, of which exactly $N^1$ are negative. Moreover, the zeros of $q_{n+1}^\alpha$ interlace the zeros of $q_n^\alpha$.

\begin{proof}[Proof of Step 4]
Since all the zeros of $P$ are real, we have $N^{\nr}=0$ and then the Step 3 implies that $q_n^\alpha$ has only real zeros, of which exactly $N^1$ are negative. Moreover, using Part (4) of Lemma \ref{nulu} we deduce that for $n\ge K$, the real zeros of $q_n^\alpha$ are always simple. Hence, the interlacing properties follows from Lemma \ref{enze}.
This proves the Part (1) of Theorem \ref{hithis}.
\end{proof}

\medskip

\noindent
\textit{Step 5.} If $N^{\nr}>0$ then there exists $n_0$ such that for $n\ge n_0$,
the polynomial $q_n^\alpha$ has exactly $n-N^{\nr}$ real (and simple) zeros, of which exactly $N^1$ are negative.

\begin{proof}[Proof of Step 5]
Using the asymptotic (\ref{eapu}), we deduce that there exists $n_0$ such that for $n\ge n_0$, $q_n$ has at least $N^{\nr}$ no real zeros. Hence we conclude from the Step 3 that $q_n^\alpha$ has exactly $n-N^{\nr}$ real zeros, and then exactly $N^1$ negative zeros.

Proceeding as in Step 3, we take $u$ big enough such that $q_n^{\alpha+u}$ has exactly $N^{\nr}$ non-real zeros. Since $q_n^\alpha$ has also exactly $N^{\nr}$ non-real zeros, we deduce from Parts (1) and (3) of Lemma \ref{nulu} that $q_n^{\alpha+i}=\Upsilon_{\alpha+i}(q_n^{\alpha+i+1})$, $0\le i\le u-1$,
has also exactly $N^{\nr}$ non-real zeros. Using Part (4) of Lemma \ref{nulu} we conclude that
the real zeros of $q_n^\alpha$ are simple.
\end{proof}

\medskip

\noindent
\textit{Step 6.} For $n$ big enough the zeros of $q_{n+1}^\alpha$ interlace the zeros of $q_n^\alpha$.

\begin{proof}[Proof of Step 6]
If we prove that there exists $n_*$ such that for $n\ge n_*$,
$q_n^\alpha$ has exactly $n-N^{\nr }$ real and simple zeros and the real zeros of the polynomial $q_{n+1}^\alpha+\lambda q_n^\alpha$ are simple for every real number $\lambda$, then the interlacing property will follow as a consequence of
Theorem \ref{obre}.

Take then $n_0$ as in the Step 5, so that $q_n^\alpha$ has exactly $n-N^{\nr}$ real (and simple) zeros, of which exactly $N^1$ are negative.

Using the notation of Remark \ref{exp} write $B_j=\gamma_j$ and define $A_j$ by
$$
P_\lambda (x)=(x+\lambda)P(x)=\sum_{j=0}^{K+1}A_jx^{K+1-j}.
$$
Then (\ref{idco}) gives
\begin{equation}\label{qgs2}
q_n^\alpha+\lambda q_{n-1}^\alpha=q_{n}^{A}=\sum_{j=0}^{K+1}A_j\La_{n-j}^\alpha.
\end{equation}
Define finally
\begin{align*}
f_n(z)&=\frac{(-1)^n\sqrt{2\pi }n^{\alpha+1/2} q_n(nz)}{(\varphi(z)e^{\frac{z}{1+\varphi(z)}})^n},\\
g_n(z)&=\frac{(-1)^n\sqrt{2\pi }n^{\alpha+1/2} q_{n-1}(nz)}{(\varphi(z)e^{\frac{z}{1+\varphi(z)}})^n},\\
f(z)&=\frac{(-1)^K\Gamma(1+\alpha)(1+\varphi(z))^\alpha}{\varphi(z)^{-1/2+K+1}z^\alpha(z^2-4z)^{1/4}}
P(-\varphi(z)).
\end{align*}

The asymptotic (\ref{asqnu}) for $q_n$ gives
\begin{equation}\label{asy2}
\lim_nf_n(z)=\varphi(z)f(z).
\end{equation}
In turns, from (\ref{qgs2}), the asymptotics (\ref{asqnu}) for $q_n^A$ and (\ref{asy2}) for $q_n$ we deduce
$$
\lim_ng_n(z)=\frac{1}{\lambda}\left(-(-\varphi(z)+\lambda)f(z)-\varphi(z)f(z)\right)=-f(z).
$$
Lemma \ref{mamon} guarantees the existence of a positive integer $n_*$ which does not depend on $\lambda$, and which can be taken $n_*\ge n_0$, such that for $n\ge n_*$
$$
f_n(z)+\lambda g_n(z)=\frac{(-1)^n\sqrt{2\pi }n^{\alpha-1/2} q_n^A(nz)}{(\varphi(z)e^{\frac{z}{1+\varphi(z)}})^n}
$$
has at least $N^{\nr}$ non-real zeros.

Hence $q_n^A$ has at least $N^{\nr}$ non-real zeros for $n\ge n_*$.  Proceeding as in the Step 5, we can prove that the real zeros of $q_n^A$ are simple. \end{proof}

The proof of the Theorem \ref{hithis} is now complete.
\end{proof}

\section{Linear combinations of Laguerre polynomials}\label{secll}
In this Section, for a positive integer $K$, we consider finite linear combinations of Laguerre polynomials of the form
\begin{equation}\label{qnla}
q_n^{\alpha;\gamma}(x)=\sum_{j=0}^K\gamma_jL_{n-j}^\alpha(x),\quad n\ge K,
\end{equation}
where we have written $q_n^{\alpha;\gamma}$ to stress the dependence on $\alpha$ and the real numbers $\gamma_j$, $1\le j\le K$.

The complete description for the zeros of $q_n^{\alpha;\gamma}$ in Theorem \ref{hithi} depends on the zeros of the polynomial
\begin{equation}\label{polpla}
P_\gamma(x)=\sum_{j=0}^K\gamma_jx^{K-j}.
\end{equation}

The proof of Theorem \ref{hithi} will need some previous results.

The following two quasi-spectral formulas for the Laguerre polynomials are one of the keys to prove Theorem \ref{hithi}
\begin{equation}\label{eapm4}
(L_n^\alpha)'(x)=-L^{\alpha+1}_{n-1}(x),\quad L_n^{\alpha-1}(x)=L_n^\alpha(x)-L_{n-1}^\alpha(x).
\end{equation}

Indeed, a straightforward use of (\ref{eapm4}) gives the following decomposition.

\begin{Lem}\label{aux} Let $\theta$ be a real zero of the polynomial $P_\gamma$. Define
then the real numbers $\tilde \gamma_{j}$, $0\le j\le K-1$, as follows
$$
\tilde \gamma_{j}=\sum_{i=0}^{j}\theta^i\gamma_{j-i}.
$$
Then
\begin{equation}\label{fesp}
q_n^{\alpha;\gamma}(x)=q_n^{\alpha-1;\tilde \gamma}(x)-(1-\theta)(q_n^{\alpha-1;\tilde \gamma})'(x).
\end{equation}
\end{Lem}

Let us note that according to Lemma \ref{lem2}, the polynomial $P_{\tilde \gamma}$ (\ref{polpla}) associated to the numbers $\tilde \gamma$ is
\begin{equation}\label{pgtg}
P_{\tilde \gamma}(x)=P_\gamma(x)/(x-\theta).
\end{equation}
(with the notation of Lemma \ref{lem2}: $A_i=\gamma_i$, $B_i=\tilde \gamma_i$).

\bigskip

As a consequence of the asymptotic (\ref{eapm5}) and (\ref{eapm6}), the sequence $q_n^{\alpha;\gamma}$ (\ref{qnla}) satisfy the following asymptotic (the proof is omitted because is the same as that of Lemma \ref{aux20}).

\begin{Lem}\label{aux2} For $z\in \CC\setminus [0,4]$, we have the following asymptotic for the polynomials $q_n$ (\ref{qnla})
\begin{equation}\label{asqn}
\lim_{n\to \infty}
\frac{(-1)^n\sqrt{2\pi n}q_n^{\alpha;\gamma}(nz)}{(\varphi(z)e^{\frac{z}{1+\varphi(z)}})^n}
=\frac{(-1)^K\varphi(z)^{1/2-K}(1+\varphi(z))^\alpha}{z^\alpha(z^2-4z)^{1/4}}
P(-\varphi(z))
\end{equation}
uniformly in compact sets of $\CC\setminus [0,4]$, where $\varphi$ is the conformal mapping of $\CC\setminus[0,4]$ onto $\{w:|w|<1\}$ defined in (\ref{defp}).
\end{Lem}

We still need the following technical result.

\begin{Lem}\label{sign} Assume $P_\gamma(1)\not =0$ and $\alpha\ge -1$ and denote
\begin{equation}\label{n0}
k_0=\left\lfloor K+1+\frac{|\alpha|}{2^{1/K}-1}\max\left\{1,\frac{\sum_{j=0}^{K-1}|\gamma_j|}{|P_\gamma(1)|}\right\}\right\rfloor.
\end{equation}
 Then for $n\ge k_0$ we have $\sign (q_n^{\alpha;\gamma}(0))=\sign(P_\gamma(1))$.
\end{Lem}

\begin{proof}
Using that $L_n^\alpha(0)=(\alpha+1)_n/n!$, we have
\begin{equation}\label{qnc}
q_n(0)=\frac{(\alpha+1)_{n-K}}{(n-K)!}\sum_{j=0}^K\left(1+\frac{\alpha}{n-j}\right)\cdots \left(1+\frac{\alpha}{n-K+1}\right)\gamma_j
\end{equation}
(to simplify the notation we do not stress the dependence on $\alpha$ and $\gamma$ in the proof of the Lemma).

Hence, we deduce
\begin{equation}\label{pce}
\lim_n\frac{(n-K)!q_n(0)}{(\alpha+1)_{n-K}}=\sum_{j=0}^K\gamma_j=P(1).
\end{equation}

Write
$$
A_n=\sum_{j=0}^K\gamma_j\prod_{i=j}^{K-1}\left(1+\frac{\alpha}{n-i}\right).
$$
(\ref{qnc}) gives that $\sign (q_n(0))=\sign(A_n)$.
We then have for $0\le j\le K-1$
$$
\prod_{i=j}^{K-1}\left(1+\frac{\alpha}{n-i}\right)=1+\sum_{l=1}^{K-j}\alpha^l\sum_{j\le i_1<\dots<i_l\le K-1}\prod_{r=1}^l\frac{1}{n-i_r}.
$$
If $n\ge K+|\alpha| s$, we have $n-i\ge |\alpha| s$, $1\le i\le K-1$, and therefore
\begin{align*}
|\epsilon_{n,j}|&=\left|\sum_{l=1}^{K-j}\alpha^l\sum_{j\le i_1<\dots<i_l\le K-1}\prod_{r=1}^l\frac{1}{n-i_r}\right|\\
&\le\frac{1}{n-K}\sum_{l=1}^{K-j}|\alpha|^l\binom{K-j-1}{l}\frac{1}{(|\alpha| s)^{l-1}}\\
&=\frac{|\alpha| s}{n-K}\sum_{l=1}^{K-j}\binom{K-j-1}{l}\frac{1}{s^{l}}=\frac{|\alpha| s}{n-K}\left[\left(1+\frac1s\right)^{K-j}-1\right].
\end{align*}
Hence, for $s=\frac{1}{2^{1/(K-j)}-1}+1$, we have
$$
|\epsilon_{n,j}|< \frac{|\alpha|}{(n-K)(2^{1/K}-1)},\quad n\ge K+\frac{|\alpha|}{2^{1/K}-1}.
$$
This gives (we take $\epsilon_{n,K}=0$)
$$
A_n=\sum_{j=0}^{K}\gamma_j(1+\epsilon_{n,j})=P(1)+\sum_{j=0}^{K-1}\gamma_j\epsilon_{n,j}.
$$
Since
$$
\left|\sum_{j=0}^K\gamma_j\epsilon_{n,j}\right|< \frac{|\alpha|}{(n-K)(2^{1/K}-1)}\sum_{j=0}^{K-1}|\gamma_j|,
$$
we deduce that for $n\ge K+\frac{|\alpha|}{2^{1/K}-1}\max\left\{1,\frac{\sum_{j=0}^{K-1}|\gamma_j|}{|P(1)|}\right\}$,
$$
\left|\sum_{j=0}^{K-1}\gamma_j\epsilon_{n,j}\right|<|P(1)|.
$$
This proves that $\sign (q_n(0))=\sign(A_n)=\sign(P(1))$.
\end{proof}

We are now ready to prove Theorem \ref{hithi}.

\begin{Rem}\label{rdlc}
We will need to consider the subset of $\RR^{K+1}$ defined by
\begin{align}\label{conA}
\Theta=\{(\gamma_0,\cdots, \gamma_K)\in \RR^{K+1}&:\gamma_0=1,\gamma_K\not=0,
\\\nonumber
&P_\gamma(z)\not=0, \mbox{ for $z\in \{z\in \CC\setminus \RR:\vert z\vert \le 1\}\cup \{1\}$}\}.
\end{align}
As it will be pointed out in the following proof, the big enough mentioned in Parts (1) and (2) of Theorem \ref{hithi} will depend uniformly on $\gamma=(\gamma_0,\cdots,\gamma_K)$ in compact sets of $\Theta$.
\end{Rem}

\begin{proof}[Proof of Theorem \ref{hithi}]

Assume first that $P(1)\not =0$. We proceed in 3 steps using induction on $K$.

\medskip

\noindent
\textit{Step 1.} The case $K=1$.

\begin{proof}[Proof of Step 1]
For $K=1$, we have that $P_\gamma$ has only one zero $\theta_1=-\gamma_1$, which it is then real. Assume first that $\alpha>0$. With the notation of Lemma \ref{aux}, we have $q_n^{\alpha-1;\tilde \gamma}(x)=L_n^{\alpha-1}(x)$, and  the identity
(\ref{fesp}) gives
$$
q_n^{\alpha;\gamma}(x)=L_n^{\alpha-1}(x)-(1-\theta_1)(L_n^{\alpha-1})'(x).
$$
Since $\alpha>0$, the Laguerre polynomial $L_n^{\alpha-1}(x)$ has $n$ positive and simple zeros. Hence, between two consecutive zeros of $L_n^{\alpha-1}(x)$, $q_n^{\alpha;\gamma}$ has an odd number of zeros.
This gives that $q_n^{\alpha;\gamma}$ has at least $n-1$ simple and positive zeros. And since the coefficients of $q_n^{\alpha;\gamma}$ are real, we deduce that it has $n$ simple and real zeros.

We consider next the case when $\theta_1<1$. If we denote $\xi_n$ for the largest zero of $L_n^{\alpha-1}(x)$, we also have
$$
\sign((L_n^{\alpha-1})'(\xi_n)=(-1)^{n},
$$
and hence $\sign (q_n^{\alpha;\gamma}(\xi_n)))=(-1)^{n+1}$. Since $\lim_{x\to +\infty}q_n^{\alpha;\gamma}(x)=(-1)^n \infty$ we deduce that $q_n^{\alpha;\gamma}$ has a zero in $(\xi_n,+\infty)$. This proves that $q_n^{\alpha;\gamma}$ has $n$ positive and simple zeros for $n\ge 1$.

If $\theta_1>1$, Lemma \ref{sign} gives that for $n_0=k_0$ given by (\ref{n0})
$$
\sign(q_n^{\alpha;\gamma}(0))=\sign(P_\gamma(1))=\sign(1-\theta_1)=-1,\quad n\ge n_0.
$$
Since $\lim_{x\to -\infty}q_n^{\alpha;\gamma}(x)=+\infty$, we deduce that for $n\ge n_0$ $q_n^{\alpha;\gamma}$ has a negative zero.  Notice that $n_0$ depends uniformly on $\gamma_1$ (in compact sets of $\Theta$).

But the result it is also true for $\alpha=0$, because $L_n^{-1}(x)=-xL_n^1(x)/n$, and we can proceed similarly as above.
\end{proof}

Assume next that the Theorem is true for $K-1$.

\medskip

\noindent
\textit{Step 2.} The case $K>1$ and $N_\gamma^{\nr }=0$.

\begin{proof}[Proof of Step 2]
If $N_\gamma^{\nr }=0$, we have that the polynomial $P_\gamma$ only has real zeros. Take then a real zero $\theta$ of $P_\gamma$. Consider the sequence $\tilde \gamma$ as defined in Lemma \ref{aux}. As a consequence of (\ref{pgtg}), we have also that
$N_{\tilde \gamma}^{\nr }=0$. Hence, the induction hypothesis gives that for $n\ge K-1$
the polynomial $q_n^{\alpha-1;\tilde \gamma}(x)$ has $n$ real and simple zeros. Hence, using the identity (\ref{fesp}), we deduce that between two consecutive zeros of $q_n^{\alpha-1;\tilde \gamma}(x)$,
$q_n^{\alpha,\gamma}$ has an odd number of zeros. And we can again conclude that actually $q_n^{\alpha,\gamma}$ has $n$ real and simple zeros.

If, in addition, $N_\gamma^{1}=0$, we have that $\theta<1$, and using the induction hypothesis, since the all the zeros of $q_n^{\alpha-1;\tilde \gamma}(x)$  are simple and positive, $q_n^{\alpha-1;\gamma}(x)$  has at least $n-1$ simple and positive zeros. As before, it is easy to prove that denoting $\xi_n$ the largest zero of $q_n^{\alpha-1;\tilde \gamma}(x)$, the polynomial $q_n^{\alpha;\gamma}(x)$ has also a zero in the interval $(\xi_n,+\infty)$. This proves that $q_n^{\alpha;\gamma}(x)$ has $n$ simple and positive zeros.

If $N_\gamma^{1}>0$,  we can take a real zero $\theta$ of the polynomial $P$ satisfying $\theta >1$.
Consider the real numbers $\tilde \gamma_j$ defined in Lemma \ref{aux}. Since $\theta>1$, we have $N_{\tilde \gamma}^1=N_{\gamma}^1-1$, $N_{\tilde \gamma}^{\nr }=N_{\gamma}^{\nr }=0$ (and hence $\alpha-1\ge K-1-N_{\tilde \gamma}^{\nr }-1$). The induction hypothesis gives that for $n$ big enough (depending uniformly on $\tilde \gamma$ in compact sets of $\Theta_{K-1}$), the polynomial $q_n^{\alpha-1;\tilde \gamma}(x)$ has $n$ real and simple zeros, of which $N_\gamma^1-1$ are negative.
Using the identity (\ref{fesp}), we deduce that between two consecutive zeros of $q_n^{\alpha-1;\tilde \gamma}(x)$, $q_n^{\alpha,\gamma}$ has an odd number of zeros.
In particular, $q_n^{\alpha,\gamma}$ has at least $N_\gamma^1-2$ negative zeros, and at least $n-N_\gamma^1$ positive zeros.

We next prove that $q_n^{\alpha,\gamma}$ has two more negative zeros, one in the interval $(-\infty,\xi_1)$ and the other in the interval $(\xi_{N^1_{\tilde\gamma}},0)$, where $\xi_1$ and $\xi_{N^1_{\tilde\gamma}}$ denote the smallest and the largest negative zero of $q_n^{\alpha-1;\tilde \gamma}(x)$, respectively.  Indeed, on the one hand, taking into account that
$$
\lim_{x\to -\infty}q_n^{\alpha,\gamma}(x)=+\infty,\quad \lim_{x\to -\infty}q_n^{\alpha-1,\tilde \gamma}(x)=+\infty,
$$
we deduce from (\ref{fesp}) that $\sign(q_n^{\alpha,\gamma}(\xi_1))=-1$, and hence we conclude that $q_n^{\alpha,\gamma}$ has one more negative zero in the interval $(-\infty,\xi_1)$. On the other hand, using again (\ref{fesp}), we have
$$
\sign(q_n^{\alpha,\gamma}(\xi_{N^1_{\tilde\gamma}}))=\sign((q_n^{\alpha-1,\tilde\gamma})'(\xi_{N^1_{\tilde\gamma}}))=(-1)^{N^1_{\tilde\gamma}}.
$$
Since for $n$ big enough
$\sign(q_n^{\alpha,\gamma}(0))=\sign(P_\gamma(1))$ (see (\ref{pce})), we deduce (taking into account that $N^1_\gamma$ is the number of real zeros greater than 1 of the monic polynomial $P_\gamma$) that $\sign(q_n^{\alpha,\gamma}(0))=\sign(P_\gamma(1))(-1)^{N^1_\gamma}=-(-1)^{N^1_{\tilde\gamma}}$. This implies that $q_n^{\alpha,\gamma}$ has one more negative zero at $(\xi_{N^1_{\tilde\gamma}},0)$. This proves that the zeros of $q_n^{\alpha-1;\tilde \gamma}(x)$ interlace the zeros of $q_n^{\alpha,\gamma}$.
\end{proof}

The interlacing properties follow now from Lemma \ref{enze}.
Hence, this proves the Part (1) of Theorem \ref{hithi}.

\medskip

\noindent
\textit{Step 3.} The case $K>1$ and $N_\gamma^{\nr }>0$.

\begin{proof}[Proof of Step 3]
We proceed by induction on $K-N_\gamma^{\nr }$.
If $N_\gamma^{\nr }=K$, as an easy consequence of Remark \ref{remw}, the asymptotic (\ref{asqn}) and the Hurwitz Theorem, we can conclude that for $n$ big enough (depending uniformly on $\gamma$ in compact sets of $\Theta$), $q_n^{\alpha;\gamma}$ has at least $N_\gamma^{\nr }=K$ non-real zeros. Theorem \ref{BeDr} implies then that $q_n^{\alpha;\gamma}$ has exactly $K$ non-real zeros and $n-K$ positive zeros.

If $K-N_\gamma^{\nr }>0$, then the polynomial $P_\gamma$ has at least one real zero $\theta$. Consider the real numbers $\tilde \gamma_j$ defined in Lemma \ref{aux}.
The identity (\ref{fesp}), together with Remark \ref{remw}, the asymptotic (\ref{asqn}) and the Hurwitz Theorem, show that for $n$ big enough (depending uniformly in compact sets of $\Theta$),
$q_n^{\alpha;\gamma}$ has at least $N_\gamma^{\nr }$ non-real zeros. If we assume that $\theta<1$, we have $N_{\tilde \gamma}^1=N_{\gamma}^1$, $N_{\tilde \gamma}^{\nr }=N_{\gamma}^{\nr }$ (and hence $\alpha-1\ge K-1-N_{\tilde \gamma}^{\nr }-1$). The induction hypothesis gives that for $n$ big enough, the polynomial $q_n^{\alpha-1;\tilde \gamma}(x)$ has $n-N_\gamma^{\nr }$ real and simple zeros, of which $N_\gamma^1$ are negative. Using the identity (\ref{fesp}), we deduce that
between two consecutive zeros of $q_n^{\alpha-1;\tilde \gamma}(x)$, $q_n^{\alpha,\gamma}$ has an odd number of zeros.
In particular, $q_n^{\alpha,\gamma}$ has at least $N_\gamma^1-1$ negative zeros, and at least $n-N_\gamma^{\nr }-N_\gamma^1-1$ positive zeros. On the one hand, denoting by $\xi_n$ the largest zero of $q_n^{\alpha-1;\tilde \gamma}(x)$, we deduce that $\sign(q_n^{\alpha,\gamma}(\xi_n))=(-1)^{n+1}$; taking into account that $\lim_{x\to +\infty}q_n^{\alpha,\gamma}(x)=(-1)^n \infty$, we conclude that $q_n^{\alpha,\gamma}$ has one more positive zero in the interval $(\xi_n,+\infty)$.
On the other hand, denoting by $\xi_{N^1_{\tilde\gamma}}$ the largest negative zero of $q_n^{\alpha-1;\tilde \gamma}(x)$, we deduce that $\sign(q_n(\xi_{N^1_{\tilde\gamma}}))=-(-1)^{N^1_\gamma}$. Since for $n$ big enough
$\sign(q_n^{\alpha,\gamma}(0))=\sign(P_\gamma(1))$ (see (\ref{pce})), we deduce (taking into account that $N^1_\gamma$ is the number of real zeros of $P_\gamma$ greater than 1) that $\sign(q_n^{\alpha,\gamma}(0))=(-1)^{N^1_\gamma}$. This implies that $q_n^{\alpha,\gamma}$ has one more negative zero at $(\xi_{N^1_{\tilde\gamma}},0)$.

The case $\theta >1$ can be proved proceeding similarly as before.
\end{proof}

The interlacing properties when $N_\gamma^{\nr }>0$ can now be proved proceeding as in the Step 6 of the proof of Theorem \ref{hithis}.

\medskip

We finally consider the case when $P(1)=0$. Assume then than $P$ has a zero at $x=1$ of multiplicity $s\le K$. Hence $P(x)=(x-1)^s\hat P(x)$, where $\hat P(x)=\sum_{j=0}^{K-s}\hat \gamma_jx^{K-s-j}$ has degree $\hat K=K-s$ and $\hat P(1)\not =0$. Using the identity (\ref{fesp}), we deduce that
$$
q_n^{\alpha;\gamma}(x)=q_n^{\alpha-s;\hat \gamma}(x).
$$
Since the polynomials $P$ and $\hat P$ have the same number of non-real zeros and the same number of real zeros largest that $1$, we deduce that $\alpha-s\ge \hat K-N_{\hat \gamma}^{\nr }-1$, and then Theorem \ref{hithi} follows from the Steps 1, 2 and 3 (because $\hat P(1)\not =0$).
\end{proof}

Since for $-2<\alpha<-1$, the Laguerre polynomials $L^\alpha_n$, $n\ge 0$, has $n$ real zeros ($n-1$ positive and one negative), proceeding as in the proof of Theorem \ref{hithi}, Steps 1 and 2, it is easy to prove the following corollary.

\begin{Cor}\label{hqj2} Assume $K-2<\alpha<K-1$ and that the polynomial $P$ (\ref{polpla}) has only real zeros. Then for $n\ge K$, $q_n^\alpha$ has $n$ real and simple zeros.
\end{Cor}

\medskip

\begin{Rem}\label{hqj}
If $\alpha$ is not an integer and $\alpha<K-2$, we next prove that for any $n\ge K$ there are real numbers $\gamma_j$, $0\le j\le K+1$, $\gamma_0=1$, $\gamma_K\not =0$, such that the polynomial $P_\gamma$ (\ref{polpla}) has only real zeros but $q_n^{\alpha;\gamma}$ (\ref{qnla}) has non-real zeros (in particular, Corollary \ref{hqj2} can not be true if $\alpha$ is not an integer and $\alpha<K-2$). Indeed, for a positive integer $K$
consider the  polynomial
$$
P(x)=(x-1)^K=\sum_{j=0}^K(-1)^j\binom{K-j}{j}x^{K-j},
$$
and the associated polynomials (\ref{qnla})
$$
q_n(x)=\sum_{j=0}^K(-1)^j\binom{K-j}{j}L^\alpha_{n-j},\quad n\ge K.
$$
Using \cite[Identity (39), p. 192]{EMOT},
we deduce $q_n(x)=L^{\alpha-K}_n(x)$. Then \cite[Theorem 6.73]{Sze} gives that if $\alpha$ is not an integer, the number of real zeros of $q_n$ is not largest than $n+\lfloor\alpha-K+1\rfloor+1$.
Hence $0\le \lfloor\alpha-K+1\rfloor+1<\alpha -K+2$, and so Corollary \ref{hqj2} can be true.

Similarly, if $\alpha<K-1$, one can deduce (by considering the positive zeros of $q_n$) that Part (1) of Theorem \ref{hithi} can not be true.
\end{Rem}

However, we next prove that Theorem \ref{hithi} is true when $\alpha$ is an integer and $P$ has only real zeros (even if $\alpha$ is a negative integer).

\begin{Theo}\label{hithz} Assume $\alpha\in \ZZ$ and that the polynomial $P$ (\ref{polpla}) has only real zeros.
\begin{enumerate}
\item If $\alpha\ge 0$, then for $n\ge K$, $q_n^\alpha$ has $n$ real and simple zeros, and the zeros of $q_{n+1}$ interlace the zeros of $q_n$. If $N_\gamma^{1}=0$ then all the zeros are positive.
\item If $\alpha\le -1$, write $s$ for the multiplicity of $x=1$ as a zero of $P$ (so that if $P(1)\not=0$ then $s=0$). Then for $n$ big enough $q_n^\alpha$ has a zero at $x=0$ of multiplicity $-\alpha+s$, and $n+\alpha-s$ real and simple zeros, and the zeros of $q_{n+1}$ interlace the zeros of $q_n$.
\end{enumerate}
In all the cases, if $N_\gamma^{1}>0$ then for $n$ big enough, $q_n^{\alpha;\gamma}$ has exactly  $N_\gamma^{1}$ negative zeros.
\end{Theo}

\begin{proof}

We first prove the Part (1) of the Theorem. We proceed by induction on $\alpha$.
Since $L_n^\alpha=\La_n^\alpha$, $n\ge 0$, for $\alpha=0$, the result follows from Theorem \ref{hithis}.

Assume next that the result is true for $\alpha\ge 0$. Hence
$$
q_n^{\alpha+1;\gamma}(x)=\sum_{j=0}^K\gamma_jL_{n-j}^{\alpha+1}(x)=-(q_{n+1}^{\alpha;\gamma}(x))',
$$
and since $q_{n+1}^{\alpha;\gamma}(x)$ has $n+1$ real and simple zeros (because of the induction hypothesis), we conclude that the polynomial $q_n^{\alpha+1;\gamma}(x)$ has
$n$ real and simple zeros. Since the polynomial $P$ does not depend on $\alpha$, using Lemma \ref{sign} it is easy to conclude that
if $N_\gamma^{1}=0$ then all the zeros are positive and if $N_\gamma^{1}>0$ then for $n$ big enough, $q_n^{\alpha;\gamma}$ has exactly  $N_\gamma^{1}$ negative zeros.

Before going to the proof of Part (2), we need to discuss the following technical question. We can assume that $P(1)\not =0$.

Write as in (\ref{n0})
$$
k_0=\left\lfloor K+1+\frac{|\alpha|}{2^{1/K}-1}\max\left\{1,\frac{\sum_{j=0}^{K-1}|\gamma_j|}{|P_\gamma(1)|}\right\}\right\rfloor.
$$
We have that $k_0$ is a function of $(\gamma_j)_{j=1}^K$, and then $k_0$ is actually a function of the zeros $(\theta_j)_{j=1}^K$ of $P$ (we are assuming that all the zeros of $P$ are real). Write then $k_{0,i}=k_0((\theta_j)_{j=1}^i)$. The proof of Theorem \ref{hithi} then shows that
the number $n_0$ in part (2) of that Theorem can be taken equal to
\begin{equation}\label{n0s}
n_0=n_0((\theta_j)_{j=1}^K)=\max\{k_{0,i}:1\le i\le K\}.
\end{equation}

We next consider the case when $\alpha\le-1$ and $P(1)\not =0$.
Using the identity
$$
L_n^{\alpha}(x)=\frac{(n+\alpha)!}{n!}(-x)^{-\alpha}L_{n+\alpha}^{-\alpha}(x), \quad n\ge -\alpha,
$$
we rewrite
\begin{align*}
q_n(x)&=\sum_{j=0}^K\gamma_jL_{n-j}^\alpha(x)=(-x)^{-\alpha}\sum_{j=0}^K\frac{(n+\alpha-j)!}{(n-j)!}\gamma_j L_{n+\alpha-j}^{-\alpha}(x)
\\&=(-x)^{-\alpha}\frac{(n+\alpha)!}{n!}\sum_{j=0}^K\frac{(n-j+1)_{j}}{(n+\alpha-j+1)_j}\gamma_j L_{n+\alpha-j}^{-\alpha}(x)
\end{align*}
It is enough to prove that for $n$ big enough the polynomial
$$
\hat q_n(x)=\sum_{j=0}^K\frac{(n-\alpha-j+1)_{j}}{(n-j+1)_j}\gamma_j L_{n-j}^{-\alpha}(x)
$$
has $n$ real and simple zeros of which $N_\gamma^1$ are negative.

We first assume that $P$ has $K$ real and simple zeros. Write
$$
P_n(x)=\sum_{j=0}^K\frac{(n-\alpha-j+1)_{j}}{(n-j+1)_j}\gamma_jx^{K-j}.
$$
Since $\lim_nP_n(x)=P(x)$ in compact set of the complex plane, we deduce that there exists $n_1$ such that for $n\ge n_1$, the polynomial $P_n$ has $K$ simple and real zeros, which we denote by $\theta_j^{[n]}$, $1\le j\le K$.
Consider the number $n_0((\theta_j^{[n]})_{j=1}^K)$ defined in (\ref{n0s}). Since also
$$
\lim_n n_0((\theta_j^{[n]})_{j=1}^K)=n_0((\theta_j)_{j=1}^K),
$$
we conclude that there exists $n_2$ such that for $n\ge n_2$,
$$
n_0((\theta_j^{[n]})_{j=1}^K)\le n_0((\theta_j)_{j=1}^K)+1=n_3.
$$
Hence, writing $n_0=\max\{n_1,n_2,n_3\}$, we have for $n\ge n_0$ that $P_n$ has $K$ simple and real zeros and
$$
n\ge n_0((\theta_j^{[n]})_{j=1}^K).
$$
The first part of this theorem then gives that the polynomials $q_n$ has $n$ real and simple zeros of which $N_\gamma^1$ are negative.

If $P$ has multiple zeros, we can proceed similarly using a continuity argument.

If $P(1)=0$, we can proceed as in the last part of the proof of Theorem \ref{hithi}.

This complete the proof of the theorem.
\end{proof}

\section{Brenke normalization of the Laguerre polynomials}\label{secl3}
In this subsection we consider the following normalization of the Laguerre polynomials.
\begin{equation}\label{lag2nx}
\Ll_n^\alpha(x)=\frac{1}{(1+\alpha)_n}L_n^\alpha(x).
\end{equation}

For real numbers $\gamma_j$, $j=0,\cdots, K$, $\gamma_0=1$ and $\gamma_K\not =0$, we define
\begin{equation}\label{polqnl}
q_n(x)=\sum_{j=0}^K\gamma_j\Ll_{n-j}^\alpha(x),\quad n\ge K,
\end{equation}
and
\begin{equation}\label{polplis}
P(x)=\sum_{j=0}^K\gamma_jx^{K-j}.
\end{equation}
Using the generating function for the generalized Laguerre polynomials $(\Ll_n^\alpha)_n$ (see \cite[vol. I, Identity (18) p. 189]{EMOT}), we find
that the polynomials $(q_n)_n$ satisfy
\begin{equation}\label{genlx}
\sum_{n=0}^\infty q_n(x)z^n=e^{z}z^KP(1/z)\pFq{0}{1}{-}{\alpha+1}{-xz}.
\end{equation}
Hence, they are Brenke polynomials. Using the terminology of \cite{drh}, $(q_n)_n$ are the Brenke polynomials generated by $e^{z}z^KP(1/z)$ and associated to $\pFq{0}{1}{-}{\alpha+1}{-z}$.

The hypergeometric function $\pFq{0}{1}{-}{\alpha+1}{-z}$ has associated the second order differential operator
\begin{equation}\label{uu2x}
\Omega_\alpha=-\frac{d}{dx}\left(\alpha I+x \frac{d}{dx}\right),
\end{equation}
so that
$\Omega_\alpha q_n=q_{n-1}$.

\medskip
In the next Lemma, we prove some properties of the operator $\Omega_\alpha$ regarding real-rootedness of polynomials.

\begin{Lem}\label{nul}
For $\alpha>-1$, the operator $\Omega_\alpha$ preserves real-rootedness (see Definition \ref{prr}). Moreover, if $p$ is a polynomial with only positive zeros then $\Omega_\alpha p$ has also only positive zeros.
If, in addition, $\alpha \ge 0$ then
\begin{enumerate}
\item $\Omega_\alpha$ is a complex zero decreasing operator, in the sense that the number of non-real zeros of $\Omega_\alpha p$ is less than or equal to the number of non-real zeros of $p$.
\item Let $p$ be a polynomial, and assume that $p$ and $\Omega_\alpha p$ has $s$ and $s-1$ real zeros, respectively. If the real zeros of $p$ are simple and $p(0)\not =0$, then the real zeros of $q$ are simple.
\end{enumerate}
\end{Lem}

\begin{proof}
According to \cite[Th. 6.2]{drh}, $\Omega_\alpha$ preserves real-rootedness. Take now a monic polynomial $p$ of degree $k$ with only positive zeros. Write $r(x)=\alpha p(x)+xp'(x)$ and $q(x)=-r'(x)$, so that $q=\Omega_\alpha p$. Assume first that the zeros of $p$ are simple and write them as $0<\xi_1<\dots<\xi_k$. It is clear that the polynomial $r$ has an odd number of zeros in each interval $(\xi_i,\xi_{i+1})$, $1\le i\le k-1$. This proves that $q$ has at least $k-2$ positive zeros, and one more zero somewhere. On the one hand, we have that $q(0)=-(\alpha+1)p'(0)=(-1)^k$. On the other hand, the leading coefficient of $q$ is equal to $-k(\alpha+k)<0$, and hence $\lim_{x\to -\infty}q(x)=(-1)^k\infty$. This shows that $q$ has to have a even number of zeros in $(-\infty,0)$ and so all the zeros of $q$ has to be positive. If $p$ has multiple zeros, we can use an argument of continuity.

We next prove Part (1) of the lemma. Consider the linear operator $T_\alpha$, acting in the linear space of polynomials, and defined by $T_\alpha(x^j)=j(\alpha+j)x^j$. Using \cite[Part 3 of Theorem 2.4]{CrCs}, we deduce that $T_\alpha$ is a complex zero decreasing operator. Since $\Omega_\alpha p=-T_\alpha (p)/x$, we have that $\Omega_\alpha$ is also a complex zero decreasing operator.

Part (2) is a consequence of Lemma \ref{nulu}. Indeed, assume first that $\alpha>0$ and consider the first order differential operator $\Upsilon_\alpha$ (\ref{uu2}). Using Part (1) of Lemma \ref{nulu}, we deduce that if $p$ has $s$ real zeros then $\Upsilon_\alpha p$ has at least $s$ real zeros. Since $\Omega _\alpha=-\alpha \frac{d}{dx} \Upsilon_\alpha$, and we assume that $\Omega_\alpha p$ has $s-1$ zeros, we conclude that $\Upsilon_\alpha p$ has exactly $r$ real zeros. Hence they are simple because of Part (4) of Lemma \ref{nulu}. If $\alpha=0$, then $\Omega _\alpha=-\frac{d}{dx} \left(x\frac{d}{dx} \right)$, and hence it is easy to conclude.
\end{proof}

\medskip

The generating function (\ref{genlx}) gives the asymptotic
\begin{equation}\label{asyl}
\lim_n\left(\frac{z}{(n+1)(\alpha+n+1)}\right)^nn!(1+\alpha)_nq_n\left(\frac{-(n+1)(\alpha+n+1)}z\right)=e^{z}z^KP(1/z),
\end{equation}
uniformly in compact sets of the complex plane (see \cite[Th. 1.1]{drh})

We are now ready to prove Theorem \ref{ultt} in the Introduction.

\begin{proof}[Proof of the Theorem \ref{ultt}]
According to \cite[Cor. 6.1]{drh}, all the zeros $q_n$, $n\ge 0$, are real if and only if $e^zz^KP(1/z)$ is a function in the Laguerre-Pólya class (see Definition \ref{def1}), i.e., if and only if all the zeros of $P$ are real.
In which case, the zeros of $q_n$  are simple and the zeros of $q_{n+1}$ interlace the zeros of $q_n$, \cite[Th. 6.3]{drh}.
Theorem \ref{BeDr} implies that at least $n-K$ of the zeros are between the zeros of the Laguerre polynomial $L_n^\alpha$, and hence in the interval $[0,4n+2|\alpha|+3]$ (see \cite[Theorem 6.31.2]{Sze}). The asymptotic (\ref{asyl}) then implies that for $n$ big enough, $q_n$ has to have $N^+$ negative zeros and $K-N^+$ positive zeros outside the interval $[0,4n+2|\alpha|+3]$. This proves the part (1).

Assume next $N^{\nr}>0$. We assume in addition that the real zeros of $P$ are simple.
Theorem \ref{BeDr} implies that for $n\ge K$, at least $n-K$ of the zeros are simple and they are between the zeros of the Laguerre polynomial $L_n^\alpha$, and hence in the interval $[0,4n+2|\alpha|+3]$.
Using Hurwitz Theorem and the asymptotic (\ref{asyl}), we can conclude that there exists $n_0$ such that for $n\ge n_0$, $q_n$ has at least $N^{\nr}$ non-real zeros.
Similarly, the asymptotic (\ref{asyl}) also implies that there exists $n_1$ such that for $n\ge n_1$, $q_n$ has to have $N^+$ negative zeros and $K-N^+$ positive zeros outside the interval $[0,4n+2|\alpha|+3]$, and they are simple.
We only have to prove that the zeros $q_{n+1}$ interlace the zeros of $q_n$.

Assuming that the real zeros of $q_n$ are simple for $n\ge n_0$, the interlacing properties can be proved proceeding as in the proof of Theorem \ref{hithis}.
Hence it is enough to prove that the real zeros of $q_n$ are simple for $n\ge n_0$.
We can assume $n_1\ge n_0$ and that for $n\ge n_0$, $q_n(0)\not=0$ (because
$$
\lim_n(n-K)!q_n(0)=\lim_n\sum_{j=0}^K\frac{\gamma_j}{(n-j)(n-j-1)\cdots (n-K+1)}=\gamma_K\not=0).
$$
Since for $n\ge n_1$, $q_n$ has exactly $N^{\nr}$ non-real zeros and $\Omega_\alpha$ is a complex zero decreasing operator, we deduce that for $n<n_1$ $q_n$ has at most $N^{\nr}$ non-real zeros. Since for $n\ge n_0$, $q_n$ has at least $N^{\nr}$ non-real zeros, we conclude that $q_n$ has exactly $N^{\nr}$ non-real zeros for $n_0\le n$ (because $q_{n-1}=\Omega_\alpha q_n$). That is, for $n\ge n_0$, the number of real zeros of $q_n$ is exactly $n-N^{\nr}$. Since for $n\ge n_1$, the real zeros of $q_n$ are simple, we deduce using Part (2) of Lemma \ref{nul} that the real zeros of $q_n$ are also simple for $n\ge n_0$ (again because $q_{n-1}=\Omega_\alpha q_n$).

If some of the zeros of $P$ are not simple, we can use an argument of continuity.

\end{proof}

\section{The case $n=K$}\label{apen}
The case $n=K$ in previous Sections is specially interesting, because it provides sufficient conditions so that an expansion of the form $\sum_{j=0}^K\tau_jL^\alpha _j$ have only real zeros. In this case, we can only consider the Laguerre polynomials, because other normalization of them can be managed by modifying the sequence $\tau_j$. We start proving Corollary \ref{coj} in the Introduction,
which provides new examples of sequences of orthogonal polynomials $(p_n)_n$ with respect to a positive measure in the real line, such that the linear
operator $T(x_n)=p_n$ preserves real-rootedness (see \cite{piot,cha,CSW}).

\begin{proof}[Proof of the Corollary \ref{coj}]
Part (1) is a consequence of Part (1) of Theorem \ref{hithi}, Corollary \ref{hqj2} and Remark \ref{hqj}.

Part (2), (3) and (4) are consequences of Part (1) of Theorems \ref{hithz}, \ref{hithis} and \ref{ultt}, respectively.
\end{proof}

\medskip

In the rest of this Section, we compare our results to those proved by Iserles, Saff and N{\o}rsett using a completely different approach (see \cite{IsSa,IsNo,INS}).

\medskip

\noindent \textit{L.1.} It is proved in \cite[Proposition 9]{IsSa} that if the polynomial $\sum_{j=0}^K\frac{\tau_j}{j!}(\alpha_j+1-x)_j$ has only positive zeros so does $\sum_{j=0}^K\tau_jL_j^{\alpha_j}(x)$.
Using our approach, we can prove this result as follows. Proceeding similarly as in the proof of Lemma \ref{joe}, it can be proved that if $\psi_i=0$, $i\ge 1$ and $\phi_i$, $1\le i\le K$, are the zeros of $Q(x)=\sum_{j=0}^K(-1)^j\gamma_{K-j}(\alpha_j+1-x)_j$, then
$$
\Be_n^{0;\phi,\psi}(x)=(-1)^K\sum_{j=0}^K \gamma_j\hat L_{K-j}^{\alpha_{K-j}}(x)
$$
(where $(\Be_n^{0;\phi,\psi})_n$ are the polynomials defined in (\ref{dhl})). Taking into account that $\hat L_n^\alpha=(-1)^nn!L_n^\alpha$ and setting $\tau_j=(-1)^jj!\gamma_{K-j}$, we have $Q(x)=\sum_{j=0}^K\frac{\tau_j}{j!}(\alpha_j+1-x)_j$ and
$$
\Be_n^{0;\phi,\psi}(x)=\sum_{j=0}^K \tau_jL_{j}^{\alpha_{j}}(x).
$$
It is now enough to use Theorem \ref{pepl}.
\bigskip

When $\alpha_i=\alpha$, $1\le i\le K$, \cite[Proposition 9]{IsSa} gives: if all the zeros of the polynomial $\sum_{j=0}^K\frac{\tau_j}{j!}(\alpha+1-x)_j$ reside in $(0,+\infty)$ then all the zeros of $\sum_{j=0}^K\tau_jL_j^{\alpha}(x)$ are positive.
We can prove something better.

\begin{Cor} Let $\alpha >-1$ and assume that the polynomial
\begin{equation}\label{uu1}
R(x)=\sum_{j=0}^K\frac{\tau_j}{j!}(\alpha+1-x)_j
\end{equation}
has only real zeros which we denote by $\zeta_i$, $1\le i\le K$, and they are arranged in decreasing order.
Write
$$
u_1=\max\{0,\lfloor -\zeta_{K-1}+1\rfloor\},\quad u_0=\max\{0,\lfloor -\zeta_{K}+1\rfloor\}.
$$
Then the Laguerre expansion
\begin{equation}\label{uu3}
\sum_{j=0}^K\tau_j (j+1)_{u_l}(\alpha+j+1)_{u_l}L^\alpha_j(x)
\end{equation}
has only real zeros for $l=1$ and only positive zeros for $l=0$.
\end{Cor}
(If $R$ has only positive zeros, then $u_0=0$ and we recover
\cite[Proposition 9]{IsSa}).

\begin{proof}
Assume we have real numbers $\gamma_j$, $0\le j\le K$, $\gamma_0=1$, $\gamma_K\not =0$, such that the polynomial
$$
Q(x)=\sum_{j=0}^K(-1)^j\gamma_j(x)_{K-j}
$$
has only real zeros. Denote by $\theta_i$, $1\le i\le K$, the real zeros of $Q$ arranged in increasing order. Write
$$
n_1=\lfloor\theta_{K-1}-\alpha\rfloor +K,\quad n_0=\lfloor\theta_{K}-\alpha\rfloor +K.
$$
We get from part (1) of Theorem \ref{ldlc} that the polynomial
$$
\sum_{j=0}^K \gamma_{K-j} \hat L^\alpha_{n-K+j}(x)
$$
has only real zeros for $n\ge n_1$, and they are positive for $n\ge n_0$. Setting $v_l=n_l-K$, $l=0,1$, and taking again into account that $\hat L_n^\alpha=(-1)^nn!L_n^\alpha$, we deduce that
the polynomial
$$
q(x)=\sum_{j=0}^K (-1)^j\gamma_{K-j}(v_l+j)! L^\alpha_{v_l+j}(x)
$$
has only real zeros for $l=1$, and only positive zeros for $l=0$.

Part (1) of Lemma \ref{nul} then gives that the polynomial
$$
\Omega_\alpha^{v_l}(q(x))=\sum_{j=0}^K (-1)^j\gamma_{K-j}(v_l+j)! (\alpha+1+j)_{v_l}L^\alpha_{j}(x)
$$
has only real zeros for $l=1$, and only positive zeros for $l=0$.

Write finally
$$
\tau_j=(-1)^{K-j}j!\gamma_{K-j}.
$$
A simple computation gives $Q(x)=R(\alpha+1-x)$, and
\begin{align*}
&\sum_{j=0}^K (-1)^j\gamma_{K-j}(v_l+j)! (\alpha+1+j)_{v_l}L^\alpha_{j}(x)\\&\hspace{2cm}=(-1)^K\sum_{j=0}^K\tau_{j}(1+j)_{v_l}(\alpha+1+j)_{v_l}L^\alpha_{j}(x),
\end{align*}
where $R$ is the polynomial defined in (\ref{uu1}).

Hence $\zeta_i=\alpha+1-\theta_i$, and so $v_l=u_l$. This proves the Corollary.
\end{proof}

\bigskip

\noindent \textit{L.2}. Iserles and Saff proved that for $\alpha>-1$, if the polynomial
\begin{equation}\label{olr}
R(x)=\sum_{j=0}^K\frac{(1+\alpha)_j\tau_j}{j!}x^j
\end{equation}
has only real zeros and they live in the interval $(-1,1)$, then the Laguerre expansion
$\sum_{j=0}^K\tau_jL^\alpha_j(x)$ has only positive zeros, see \cite[(ii), p. 560]{IsSa}. The case $n=K$ in Part (1) of Theorem \ref{hithis} prove something much better.
Indeed, if $R(x)=\sum_{j=0}^K\frac{(1+\alpha)_j\tau_j}{j!}x^j$ has only real zeros (no matter if they live in $(-1,1)$ or not), Theorem \ref{hithis} establishes that all the zeros of $\sum_{j=0}^K\tau_jL^\alpha_j(x)$ are also real and simple. Moreover, those zeros are positive if the zeros of $P$ live in $(-\infty, 1)$.

\bigskip

\noindent \textit{L.3}. The case $n=K$ in Part (1) of Theorem \ref{hithi} and Corollary (\ref{hqj2}) gives the following.
Let $\alpha \ge K-2$. If the polynomial
\begin{equation}\label{olp}
P(x)=\sum_{j=0}^K\tau_jx^j
\end{equation}
has only real zeros, then the Laguerre expansion
$\sum_{j=0}^K\tau_jL^\alpha_j(x)$ has only real zeros. Moreover, if the zeros of $P$ are less than 1 and $\alpha\ge K-1$, then $\sum_{j=0}^K\tau_jL^\alpha_j(x)$ has only real and positive zeros. This result seems to have gone unnoticed by Iserles and Saff.

When $\alpha\in \NN$, there is a connection between (L.2) and (L.3). Indeed, consider the operator $T$ acting in the linear space of polynomials as $T(x^n)=(1+\alpha)_nx^n/n!$.
Then, if $R$ and $P$ are the polynomials (\ref{olr}) and (\ref{olp}), respectively, then $T(P)=R$.
In general, if $p$ is a real-rooted polynomial $T(p)$ does not need to be real-rooted. Indeed, this is equivalent to the sequence $((1+\alpha)_n/n!)_n$ being a multiplier (see Definition \ref{mul}), which in turn it is equivalent to (see Theorem \ref{PS2}))
$$
\sum_{j=0}^\infty \frac{(1+\alpha)_j}{j!j!}x^j=\pFq{1}{1}{1+\alpha}{1}{x}\in \lp I.
$$
However, according to \cite{KiKim} $\pFq{1}{1}{1+\alpha}{1}{x}\in \lp I$ if and only if $\alpha=m$, $m\in \NN$. In fact,
if $\alpha =m\in \NN$, we have
$$
\pFq{1}{1}{1+m}{1}{x}=m!e^xL_m(-x)\in \lp I.
$$
Hence, for $\alpha=m\in \NN$, (L.2) implies (L.3).

\bigskip

\noindent\textit{L.4}. Let $\alpha >-1$. If the polynomial $P(x)=\sum_{j=0}^K(1+\alpha)_j\tau_j x^j$ has only real zeros, then the Laguerre expansion
$\sum_{j=0}^K \tau_jL^\alpha_j(x)$ has only real zeros.

Since the sequence $(1/(1+\alpha)_n)_n$ is a multiplier (see Definition (\ref{mul})), we have that for $\alpha>K-1$, \textit{L.3} implies \textit{L.4}.

\bigskip

\section{Asymptotic for the zeros of the polynomials $q_n$}\label{bez}

In the last Section of this paper, and for the sake of completeness, we establish the asymptotic behaviour of the zeros of polynomials $(q_n)_n$ studied in the last Sections.

This asymptotic has already been established in \cite[Corollary 4.5]{Dur0} for the polynomials $q_n$ (\ref{qnml}) studied in Section \ref{ssec4} and associated to the monic Laguerre polynomials.

In the next Corollary we establish the asymptotic behaviour for the zeros of the polynomials $(q_n)_n$ (\ref{qnu}) studied in Section \ref{secl33} and associated to the normalization of the Laguerre polynomials taking the value $1$ at zero (\ref{lag2n}).

\begin{Cor}\label{pcoru} Under the hypothesis of Theorem \ref{hithis},
for $n$ big enough write $\xi_j(n)$, $1\le j\le n-N^{\nr }$, for the real zeros of the polynomial $q_n$ (\ref{qnu}) arranged in increasing order, and $\xi^{\nr }_j(n)$, $j=1,\cdots, N^{\nr }$, for the non-real zeros
arranged in increasing lexicographic order. Write also $\theta_j$, $1\le j\le K-N_\gamma^{\nr }$, for the real zeros of the polynomial $P$ (\ref{dPl2}) arranged in decreasing order, and $\theta^{\nr }_j(n)$, $j=1,\cdots, N^{\nr }$, for the non-real zeros arranged in decreasing lexicographic order. Write finally, $N^{-1}$ for the number of real zeros of $P$ smaller than $-1$. If we assume in addition that $P(1)\not =0$, then:
\begin{enumerate}
\item Asymptotic for the leftmost positive zeros: for $N^1+1\le i\le n-N^{\nr }-N^{-1}$,
\begin{equation}\label{aql1}
\lim_n n\xi_{i}(n)=j_{i-N^1,\alpha}^2/4,
\end{equation}
where $j_{i,\alpha}$ denotes the $i$-th positive zero of the Bessel function $J_\alpha$.
\item Asymptotic for the negative zeros:
\begin{equation}\label{aql2}
\lim_n \frac{\xi_{j}(n)}{n}=-\theta_j-1/\theta_j+2,\quad 1\le j\le N^{1}.
\end{equation}
\item Asymptotic for the rightmost positive zeros:
\begin{equation}\label{aql3}
\lim_n \frac{\xi_{j}(n)}{n}=-\theta_{j-n+K}-1/\theta_{j-n+K}+2,\quad n-N^{\nr }-N^{-1}+1\le j\le n-N^{\nr }.
\end{equation}
\item Asymptotic for the non-real zeros:
$$
\lim_n \frac{\xi^{\nr }_{j}(n)}{n+1}=-\theta_j^{\nr }-1/\theta_j^{\nr }+2,\quad 1\le j\le N^{\nr }.
$$
\item For a bounded continuous function $f$ in $\RR$, we have
$$
\lim_n\frac{1}{n}\sum_{j=1}^nf\left(\frac{\xi_j(n)}{n}\right)=\frac{1}{2\pi}\int_{0}^4 f(x)\sqrt{\frac{4-x}x}dx,
$$
\end{enumerate}

\end{Cor}

\begin{proof}
Using the Mehler-Heine type formula for the Laguerre polynomials
\begin{equation}\label{lel}
\lim_n\frac{L_n^\alpha(z/(n+j))}{n^\alpha}=z^{-\alpha/2}J_\alpha(2\sqrt z)=\frac{1}{\Gamma (\alpha+1)}\pFq{0}{1}{-}{\alpha+1}{-z}
\end{equation}
uniformly in compact set of the complex plane,
it is easy to deduce the following asymptotic for the polynomials
$q_n$ (because of (\ref{gan2}))
$$
\lim_nq_n(x/n)=P_\gamma(1)\Gamma(1+\alpha)z^{-\alpha/2}J_\alpha(2\sqrt z)=P_\gamma(1)\pFq{0}{1}{}{\alpha+1}{-x}
$$
This together with Theorem \ref{hithis} and the Hurwitz's Theorem proves the asymptotic (\ref{aql1}). In order to prove the other asymptotics, we use the Theorem \ref{BeDr}. This Theorem proves that $n-K$ of the real zeros of each $q_n$ are between the zeros of the Laguerre polynomial $L_n^\alpha$. And so they are positive and satisfy $\xi_j(n)<4n+2|\alpha|+3$. Since for these zeros
$$
\frac{\xi_j(n)}{n}\le 4+\frac{2|\alpha|+3}{n}\to 4,
$$
and $x+1/x>2$, for $x>1$, the Parts (2), (3) and (4) follow as a consequence of Remark \ref{remw} and the asymptotic (\ref{asqnu}).
The Part (5) also follows from Theorem \ref{BeDr} using the following well-known weak scaling limit
for the counting measure of the zeros of the Laguerre polynomials (\cite[Theorem 1]{Gaw})
\begin{equation}\label{wsl2}
\lim_n\frac{1}{n}\sum_{j=1}^nf\left(\frac{\zeta_j^\alpha(n)}{n}\right)=\frac{1}{2\pi}\int_{0}^4 f(x)\sqrt{\frac{4-x}x}dx,
\end{equation}
where  $\zeta_j^\alpha(n)$ denote the zeros of the Laguerre polynomial $L^\alpha_n$ arranged in increasing order.
\end{proof}

\medskip

The zeros of the polynomials $(q_n)_n$ (\ref{qnla}), studied in Section \ref{secll} and associated to the standard Laguerre polynomials, behave as those of the polynomials $(q_n)_n$ (\ref{qnu}) as displayed in the previous Corollary.
This is because the polynomials $q_n$ (\ref{qnla}) also satisfy both, a Mehler-Heine type formula:
$$
\lim_n\frac{1}{n^\alpha}q_n(x/n)=P_\gamma(1)z^{-\alpha/2}J_\alpha(2\sqrt z)=\frac{P_\gamma(1)}{\Gamma(1+\alpha)}\pFq{0}{1}{}{\alpha+1}{-x},
$$
and the asymptotic (\ref{asqn}).

\medskip

Finally, the asymptotic of the zeros of the polynomials $(q_n)_n$ (\ref{polqnl}), studied in Section \ref{secl3} and associated to the Brenke normalization of the Laguerre polynomials (\ref{lag2nx}), is as follows.

\begin{Cor}\label{ultc} Under the hypothesis of Theorem \ref{ultt},
for $n$ big enough write $\xi_j(n)$, $1\le j\le n-N^{\nr }$, for the real zeros of the polynomial $q_n$ (\ref{polqnl}) arranged in increasing order, and $\xi^{\nr }_j(n)$, $j=1,\cdots, N^{\nr }$, for the non-real zeros
arranged in increasing lexicographic order. Write also $\theta_j$, $1\le j\le K-N^{\nr }$, for the real zeros of the polynomial $P$ (\ref{polplis}) arranged in decreasing order, and $\theta^{\nr }_j(n)$, $j=1,\cdots, N^{\nr }$, for the non-real zeros arranged in decreasing lexicographic order. Write finally, $N^{+}$ for the number of positive zeros of $P$. Then:
\begin{enumerate}
\item Asymptotic for the leftmost positive zeros: for $N^++1\le i\le n-K+N^{+}$,
\begin{equation}\label{aql1x}
\lim_n n\xi_{i}(n)=j_{i-N^+,\alpha}^2/4,
\end{equation}
where $j_{i,\alpha}$ denotes the $i$-th positive zero of the Bessel function $J_\alpha$.
\item Asymptotic for the negative zeros:
\begin{equation}\label{aql2x}
\lim_n \frac{\xi_{j}(n)}{n^2}=-\theta_j,\quad 1\le j\le N^{+}.
\end{equation}
\item Asymptotic for the rightmost positive zeros:
\begin{equation}\label{aql3x}
\lim_n \frac{\xi_{j}(n)}{n^2}=-\theta_{j-n+K},\quad n-K+N^{+}+1\le j\le n-N^{\nr }.
\end{equation}
\item Asymptotic for the non-real zeros:
$$
\lim_n \frac{\xi^{\nr }_{j}(n)}{n^2}=-\theta_j^{\nr },\quad 1\le j\le N^{\nr }.
$$
\item For a bounded continuous function $f$ in $\RR$, we have
$$
\lim_n\frac{1}{n}\sum_{j=1}^{n-N^{\nr}}f\left(\frac{\xi_j(n)}{n}\right)=\frac{1}{2\pi}\int_{0}^4 f(x)\sqrt{\frac{4-x}x}dx,
$$
\end{enumerate}

\end{Cor}

We omit the proof because is similar to that of Corollary \ref{pcoru}.



\end{document}